\documentclass[MSNbibl,number,citesort,seceqn,dvips]{arxbj}
\usepackage{upgreek}
\usepackage{algorithm}
\usepackage{algorithmic}
\usepackage{graphicx}

\aid{0}
\volume{20}
\issue{3}
\pubyear{2014}
\firstpage{1126}
\lastpage{1164}
\doi{10.3150/13-BEJ517} 

\makeatletter

\newcommand{\angler}{\rangle}
\newcommand{\anglel}{\langle}
\newcommand{\rrvert}{\vert}
\newcommand{\llvert}{\vert}
\newcommand{\iint}{\int\!\!\!\int}
\newtheorem{thrm}{Theorem}[section]
\newremark{remk}{Remark}[section]

\newtheorem{lmma}{Lemma}[section]
\newtheorem{prop}{Proposition}[section]
\newremark{example}{Example}

\def\calF{\mathcal{F}}

\def\calL{\mathcal{L}}

\def\calI{\mathcal{I}}

\def\calL{\mathcal{L}}
\def\calK{\mathcal{K}}
\def\calN{\mathcal{N}}

\def\calT{\mathcal{T}}

\def\calR{\mathcal{R}}

\def\calX{\mathcal{X}}
\def\bbr{\mathbb{R}}
\def\bbe{\mathbb{E}} 
\def\bbp{\mathbb{P}}

\newcommand{\eqref}[1]{(\ref{#1})}
\makeatother

\begin{document}
\begin{frontmatter}

\title{Small-time asymptotics of stopped L\'evy bridges and simulation
schemes with controlled bias}
\runtitle{Small-time asymptotics of stopped L\'evy bridges}

\begin{aug}
\author[a]{\inits{J.E.}\fnms{Jos\'e E.}~\snm{Figueroa-L\'opez}\thanksref{a}\ead[label=e1]{figueroa@purdue.edu}}
\and
\author[b]{\inits{P.}\fnms{Peter}~\snm{Tankov}\thanksref{b}\ead[label=e2]{tankov@math.univ-paris-diderot.fr}}
\runauthor{J.E. Figueroa-L\'opez and P. Tankov} 
\address[a]{Department of Statistics, Purdue University, 250 N.
University Street, West Lafayette, IN 47907, USA.
\printead{e1}}
\address[b]{Laboratoire de Probabilit\'es et Mod\`eles Al\'eatoires,
Universit\'e Paris-Diderot (Paris 7), Case 7012,
75205 Paris Cedex 13, France.
\printead{e2}}
\end{aug}

\received{\smonth{3} \syear{2012}}
\revised{\smonth{1} \syear{2013}}

%
\begin{abstract}
We characterize the small-time asymptotic behavior of the exit
probability of a L\'evy process out of a two-sided interval {and} of
the law of its overshoot, conditionally on the terminal value of the
process. The asymptotic expansions are given in the form of a first-order term and a precise computable error bound. As an important
application of these formulas, we develop a novel adaptive
discretization scheme for the Monte Carlo computation of functionals of
killed L\'evy processes with controlled bias. The considered
functionals appear in several domains of mathematical finance (e.g.,
structural credit risk models, pricing of barrier options, and
contingent convertible bonds) as well as in natural sciences. The
proposed algorithm works by adding discretization points sampled from
the L\'evy bridge density to the skeleton of the process until the
overall error for a given trajectory becomes smaller than the maximum
tolerance given by the user.
\end{abstract}

%

\end{frontmatter}

\section{Introduction}\label{sec1}

Small-time asymptotics for the distributions of L\'evy processes and
related Markov processes have a long history going back to the
seminal work of L\'eandre \cite{Leandre}, who obtained the leading
order term of the transition density of a Markov process
solving a stochastic differential equation with jumps.
In the case of a L\'evy process, the main result of L\'eandre \cite{Leandre}
reads
%
%
\begin{equation}
\label{Eq:AsymDstyPureLevy} \lim_{t\to 0}\frac{1}{t}
f_{t}(x)=s(x)\qquad (x\neq 0),
\end{equation}
where $f_{t}(x):=\frac{\mathrm{d}}{\mathrm{d}x} \bbp(X_{t}\leq x)$ is the marginal
density of the L\'evy process $X$ and $s$ is the L\'evy density of
$X$, {whose existence and smoothness need to be assumed}. L\'eandre's
approach was to
consider separately the small jumps (say, those with sizes smaller
than an $\varepsilon>0$) and the large jumps of the underlying L\'evy
process, and to condition on the number of large jumps by time~$t$.
A similar approach has been applied during the last decade to obtain
high-order asymptotic expansions for the transition distributions and
densities of L\'evy processes and {other Markov processes with jumps} (see
R\"uschendorf and
Woerner~\cite{ruschendorf02}, Figueroa-L\'opez, Gong and
Houdr\'e~\cite{FigGongHou:2011}, Figueroa-L\'opez and
Houdr\'e~\cite{FigHou:2009}, and
Figueroa-L\'opez and Ouyang \cite{FLO11}). These small-time asymptotic results have found a wide
{scope} of applications ranging from estimation methods based on
high-frequency sampling observations of the process (see, e.g., %
Figueroa-L\'opez \cite{Fig:2008a},
Comte and Genon-Catalot \cite{CmtGnn:2011}, Rosenbaum and
Tankov \cite{rosenbaum.tankov.10}, {and references
therein}) to asymptotic results for option prices and Black--Scholes
volatilities in short-time (cf. Tankov~\cite{Tnkv10}, Figueroa-L\'{o}pez and Forde \cite{FLF11}, Figueroa-L\'opez, Gong and
Houdr\'e \cite
{FigGongHou:2011}).

In the present paper, we adopt Leandre's approach to study the
asymptotic behavior of the generalized moments of the L\'evy process
stopped at the time it exits a two-sided interval $(a,b)$,
conditionally on the terminal value of the process. Specifically, for a
L\'evy process $(X_{t})_{t\geq 0}$ with L\'evy density $s$ that is
smooth outside any neighborhood of the origin and for a bounded
Lipschitz function $\varphi$, we prove that
%
%
\begin{equation}
\label{Eq:IntroMainAsympRes} \bbe \bigl(\varphi(X_{\tau}){
\mathbf{1}}_{\tau\leq t}\vert X_{t}=y \bigr)=\frac{t}{2} \int
_{(a,b)^{c}}\varphi(v) \frac
{s(v)s(y-v)}{s(y)} \,\mathrm{d}v+\mathrm{o}(t) \qquad(t\to 0)
\end{equation}
for any $y\in(a,b)\setminus\{0\}$, where $\tau:=\inf\{u\geq 0\dvt X_{u}\notin(a,b)\}$ with $-\infty\leq
a<0<b\leq\infty$.
In the case $\varphi\equiv 1$, (\ref{Eq:IntroMainAsympRes}) can be
written as follows:
%
%
\begin{equation}
\label{Eq:AsympExitProb} \bbp \bigl(\exists u\in[0,t]\dvt X_{u}
\notin(a,b)\vert X_{t}=y \bigr)=\frac{t}{2} \int
_{(a,b)^{c}}\frac{s(v)s(y-v)}{s(y)} \,\mathrm{d}v+\mathrm{o}(t) \qquad (t\to 0)
\end{equation}
for {$y\in(a,b)\setminus\{0\}$}. As in the case of the small-time
asymptotics for the
{marginal} distributions of the process, the main intuition can be
drawn from
considering the pure-jump case with finite jump {activity.
Intuitively,} formulas \eqref{Eq:IntroMainAsympRes}--\eqref
{Eq:AsympExitProb} tell us that if, within a small time period, a L\'
evy process goes out of the interval $(a,b)$ and then comes back to the
point $y\in(a,b)$, this essentially happens with two large jumps:
the first {jump takes} the process out of $(a,b)$, {while} the second
{jump brings} it back to $y$.

Our study of the short-time behavior of
(\ref{Eq:IntroMainAsympRes}) and (\ref{Eq:AsympExitProb})
is motivated {by applications in the} Monte Carlo {evaluation} of
functionals of the form 
%
%
\begin{eqnarray}
\bbe\bigl[{F}(X_T){\mathbf{1}}_{\tau> T}\bigr],\qquad \tau= \inf\bigl
\{t\geq0\dvt X_t \notin (a,b)\bigr\}. \label{killedlevyintro}
\end{eqnarray}
In financial mathematics, such functionals arise in structural credit
risk models based on L\'evy processes (Fang \textit{et~al.} \cite{fang.al.10}) and in the
pricing of barrier options (cf. Kou and Wang \cite{KouWang2000},
Boyarchenko and
Levendorskii \cite{BoyLev2002}), which is one of the most popular classes of exotic
options. Very recently, a renewed interest to these problems has
emerged in relation to the {so-called} contingent convertible bonds,
where the conversion is triggered by a passage across a level and
which exhibit a high sensitivity to jump risk
(Corcuera \textit{et~al.} \cite{corcuera2011efficient}).
In natural sciences, L\'evy processes (under the name of L\'evy
flights) are used as models for certain diffusion-like phenomena
in physics and chemistry (so-called anomalous or super-diffusion)
(Metzler and Klafter~\cite{metzler2000random},
Shlesinger, Zaslavsky and
Frisch \cite{shlesinger.al.95},
Barthelemy, Bertolotti and
Wiersma \cite{barthelemy2008levy})
as well as to describe movement patterns of foraging animals
(Viswanathan
\textit{et~al.} \cite{viswanathan1996levy},
Benhamou \cite{benhamou2007many}), and there is considerable
interest toward the study of
L\'evy flights in bounded domains and related first passage problems
giving rise to functionals of type (\ref{killedlevyintro}) (Chechkin \textit{et~al.} \cite
{chechkin2005barrier},
Buldyrev
\textit{et~al.} \cite{buldyrev2001properties},
Garbaczewski and
Stephanovich \cite{garbaczewski2009levy}).
In all these settings, closed-form expressions are rarely
available {and Monte} Carlo is often the method of choice.

The simplest procedure to evaluate the functional
\eqref{killedlevyintro} by Monte Carlo consists in simulating the
process $(X_{t})_{t\geq 0}$ at evenly spaced times $t_{k}^{n}:=k
h_{n}$, with $h_{n}:=T/n$ {and $k=0,\ldots,n$}, over the {interval
$[0,T]$,} and approximating the exit time {$\tau$} by
\[
\tilde{\tau}_{n}:=\inf\bigl\{t_{k}^{n}\dvt X_{t_{k}^{n}}\notin({a},{b})\bigr\}.
\]
This simple method introduces two types of errors: the statistical
error and the discretization error. The latter is known to be quite
significant (cf. Baldi \cite{Baldi} and Example \ref{ex2} in Section~\ref{Sec:NumIllstr} below); Metwally and Atiya \cite{MetAti2002} reports errors of up to
$10\%$ in the context of barrier options for a time discretization of
one point per day.

In the context of continuous diffusions, short-time asymptotics have
been successfully employed to alleviate the bias due to
the discretization error. One of the earliest procedures of this type,
due to Baldi \cite{Baldi}, is based on an approximation of the probability, $p(x,y,t)$,
that the process $X$ has gone out of a domain $(a,b)$
{during the small time interval $[s,s+t]$ conditioning
on $X_{s}=x$ and $X_{s+t}=y$; that is,}
%
%
\begin{equation}
\label{Eq:GenDfnExitProba} p(x,y,t):=\bbp \bigl(\exists u\in[s,s+t]\dvt X_{u}\notin (a,b)\vert X_{s}=x,X_{s+t}=y
\bigr).
\end{equation}
Given such an
approximation $\tilde{p}(x,y,t)$ of the functional $p(x,y,t)$, the
procedure simulates iteratively $X_{t_{k+1}^{n}}$ at each step
$k=0,\ldots,n-1$, and if $X_{t_{k+1}^{n}}\in(a,b)$, it proceeds to
kill the process with probability
$\tilde{p}(X_{t_{k}^{n}},X_{t_{k+1}^{n}},h_{n})$ and choose
$t_{k+1}^{n}=(k+1)h_{n}$ as an approximation of the exit time
$\tau$. A similar idea was used in Moon \cite{Moon2008} to price barrier
options with payoff $\varphi(S_{\tau},\tau)$ by Monte Carlo.

In the context of L\'evy processes, an attempt to apply a similar
methodology has been made in Webber \cite{Webber},
Ribeiro and
Webber \cite{ribeiro2006correcting}. The
authors {remarked} that the discretization bias can be {reduced} by
using the identity
%
%
\begin{equation}
\label{Eq:KyDcmpDetermTimes} \bbe \bigl({F} (X_T) {
\mathbf{1}}_{\{\tau<T\}} \bigr)=\bbe \Biggl({F} (X_T) \Biggl(1-\prod
_{k=0}^{n-1} \bigl\{ 1-p(X_{t_{k}^{n}},X_{t_{k+1}^{n}},h_{n})
\bigr\} \Biggr) \Biggr)
\end{equation}
and replacing the exact exit probability $p(x,y,t)$ with a suitable
small-time approximation
$\tilde{p}(x,y,t)$. However, these papers propose no general formula
for $\tilde p(x,y,t)$ and, as shown in Becker \cite{becker.10}, the Monte
Carlo method proposed in Webber \cite{Webber},
Ribeiro and
Webber \cite{ribeiro2006correcting} could lead
to a large discretization bias. On the other hand, in the specific case
of the parametric
{variance} gamma model, there exist discretization algorithms
{(cf.~Avramidis and
L'Ecuyer \cite{avramidis2006efficient})} allowing to simulate the running minimum
and maximum with error bounds. Let us also remark the {recent} work of
Kuznetsov \textit{et~al.} \cite{Kuznetsovetal2011} where a {method for the joint simulation of
the running maximum and the position} of a L\'evy process is introduced
based on the Wiener--Hopf decomposition of the process.

Our short-time asymptotic result (\ref{Eq:AsympExitProb}) provides an
approximation of the exit probability {(\ref{Eq:GenDfnExitProba})} via
the formula
%
%
\begin{equation}
\label{Eq:DfnApproxExitProba0} {\tilde{p}(x,y,t):=\frac{t}{2} \int
_{(a-x,b-x)^{c}}\frac
{s(v)s(y-x-v)}{s(y-x)} \,\mathrm{d}v=\frac{t}{2} \int
_{(a,b)^{c}}\frac
{s(u-x)s(y-u)}{s(y-x)} \,\mathrm{d}u}
\end{equation}
for $x\neq y$, which is valid under mild regularity conditions on the
L\'evy process $X$ (see Section~\ref{asymp} for details). {The first-order approximation (\ref{Eq:DfnApproxExitProba0}), together with an
appropriate error bound for it, enable} us to develop a general
adaptive Monte
Carlo method for evaluating the functional \eqref{killedlevyintro} with
a given precision. Given a target error level {$\gamma$}, the idea is to
generate a ``random skeleton'' $\{(T_{k},X_{T_{k}})\}_{k=1}^{N}$ of
the process $X$ such that the error in each subinterval
$[T_{k},T_{k+1}]$, that is,
%
%
\begin{equation}
\label{Eq:IniErrorDfn} e:=p(X_{T_{k}},X_{T_{k+1}},T_{k+1}-T_{k})-
\tilde {p}(X_{T_{k}},X_{T_{k+1}},T_{k+1}-T_{k}),
\end{equation}
satisfies $|e|\leq\frac{T_{k+1}-T_k}{T} {\gamma}$. The
functional \eqref{killedlevyintro} is then approximated as {follows:}
%
%
\begin{equation}
\label{Eq:KyDcmpFrm} \bbe\bigl[{F} (X_T){\mathbf{1}}_{\tau> T}
\bigr] \approx\bbe \Biggl({F} (X_T) \prod
_{k=0}^{N-1} \bigl\{1-\tilde {p}(X_{T_{k}},X_{T_{k+1}},T_{k+1}-T_k)
\bigr\} \Biggr),
\end{equation}
and it is shown that the total bias of this computation will be less
then {$\gamma$}. As a result of this adaptiveness, the algorithm generates
{more frequent} points when the process $X$ is close to the boundary,
and takes
large {time} steps (thus saving computational time) when the process is far
from the boundary. {Let us remark that, unlike the formula (\ref
{Eq:KyDcmpDetermTimes}), where the sampling times $\{t_{k}^{n}\}$ are
deterministic and fixed, the decomposition (\ref{Eq:KyDcmpFrm}) for
random skeletons $\calX:=\{(T_{k},X_{T_{k}})\}_{k=0}^{N}$ requires
precise (and also novel to the best of our knowledge) conditions under
which this formula holds (see Section~\ref{adaptive.sec} for the details).}

The {proposed adaptive} algorithm works as follows. First, the endpoint
$X_T$ is generated
and added to the skeleton. Next, if {the error (\ref{Eq:IniErrorDfn})}
is too large for a
given subinterval {$[T_{k},T_{k+1}]$}, the {procedure} splits the
interval into two and
generates the midpoint $X_{\bar{T}_{k}}$ with
$\bar{T}_{k}:=(T_{k}+T_{k+1})/2$ from the bridge distribution. This is
repeated {iteratively until the desired error bound is satisfied for
every subinterval {$[T_k,T_{k+1}]$} of the sampling times
$0=T_{0}<\cdots<T_{N}=T$}. Such retrospective {sampling} (starting from the
endpoint) has a number of advantages over the classical uniform
discretization, especially in the context of rare event simulation,
where it enables one to easily implement variance reduction by
importance sampling. Indeed,
the process can be directed to the region of interest by modifying the
distribution of the terminal value, while keeping unchanged the rest
of the algorithm. On the
other hand, this method requires fast simulation from the bridge
distribution of $X_{t/2}$ conditioned to $X_{t}=y$. {To this end}, as {another}
contribution {of particular interest on its own}, we {also} propose
a new method to simulate from this {L\'evy} bridge distribution based
on the classical rejection method.

{As previously explained, in order to} implement the above adaptive
algorithm, {precise computable bounds for the approximation errors in
(\ref{Eq:IntroMainAsympRes})--(\ref{Eq:AsympExitProb}) are also needed.
We obtain such bounds by developing explicit inequalities for the tail
probabilities and transition densities of a L\'evy process whose L\'evy
density has a small compact support. This type of concentration
inequalities in turn allows us to estimate} the different components of
the error, which, as explained above, {originate} from conditioning the
desired functional on the number of big jumps by time $t$ {(see Section~\ref{estgen} for the details)}.
The resulting error {bounds are} given in terms of the Lipschitz and
$L_{\infty}$ norms of $\varphi$ as well as several computable
quantities related to the L\'evy density $s$ such as
$\sup_{|x|\geq \varepsilon} s(x)$, $\sup_{|x|\geq \varepsilon}
|s'(x)|$, $\int_{|x|\geq \varepsilon} s(x) \,\mathrm{d}x$, and
$\int_{|x|\leq\varepsilon} x^{2}s(x)\,\mathrm{d}x$.

{Let us also remark that} {an} adaptive simulation method similar to
the one introduced in the
present paper was proposed in Dzougoutov \textit{et~al.} \cite{Dzouetal2005} to compute a
functional of the form $\bbe\varphi(X_{\tau},\tau)$ for a homogeneous
diffusion process $X$ {without jumps}. Adaptive numerical methods for
finding weak approximation of diffusions without jumps and with finite
intensity jumps (but with the adaptiveness only concerning the
diffusion part) have also been proposed in Szepessy, Tempone and
Zouraris \cite{Szepessy2001} and Mordecki \textit{et~al.} \cite
{Mordecki}, respectively. As in our paper, the idea therein is to
sample from inside of a subinterval $[t_{k}^{n},t_{k+1}^{n}]$ whenever
the approximation error in that subinterval has not reached a desired
low level, specified by the user.

The paper is organized as follows. In Section~\ref{asymp}, we obtain
the leading term of the functional $\bbe (\varphi(X_{\tau
})\mathbf{1}_{\tau\leq t}\vert  X_{t}=y )$ when $t\to 0$. The
{explicit} estimate of the approximation error is given in Section~\ref{estgen}.
The development of the adaptive discretization schemes for the Monte Carlo
computation of the functional {$\bbe[{F} (X_T){\mathbf{1}}_{\tau> T}]$}
as well as the algorithm to simulate random observations from the L\'
evy bridge distribution are given in
Section~\ref{adaptive.sec}. Our methods are illustrated numerically in
Section~\ref{Sec:NumIllstr} {for} Cauchy process. Finally, the proofs
of the technical results are deferred to the \hyperref[app]{Appendix}.

\section{Small-time asymptotics for L\'evy bridges}\label{asymp}
Let $X$ be a real-valued L\'evy process on a probability space $(\Omega
,\calF,\bbp)$ with L\'evy triplet
$(\sigma^{2},\nu,\mu)$ with respect to truncation function
$h(x) = \mathbf{1}_{|x|\leq1}$. Throughout, $(\calF_{t})_{t\geq 0}$ denotes
the natural filtration generated
by the process $X$ and augmented by the null sets of $\calF$ so
that it satisfies the usual conditions (see, e.g., Chapter I.4 in
Protter \cite{Protter}). The following standing assumptions are imposed
throughout the paper:
\begin{itemize}
\item The L\'evy measure $\nu$ admits a continuously
differentiable density $s\dvtx \bbr\setminus\{0\}\to(0,\infty)$, with
respect to the Lebesgue measure (hereafter denoted by $\calL$), which
satisfies, for any $\varepsilon>0$,
%
%
\begin{equation}
\label{Cond1} \sup_{|x|\geq \varepsilon} s(x)<\infty,\qquad \sup
_{|x|\geq
\varepsilon} \bigl|s'(x)\bigr|<\infty.
\end{equation}
\item The distribution of $X_{t}$ admits a density $f_{t}$ for all
$t> 0$. Since $\nu$ {is already assumed to} admit a density, for this
assumption {to hold,} it suffices to additionally
require that $\nu(\bbr)=\infty$ or $\sigma>0$
(see Theorem 27.7 in Sato~\cite{sato}).
\item The density of $X_t$ satisfies $f_{t}(x)>0$ for all $x\in\bbr$
and $t>0$ (see Theorem 24.10 in Sato~\cite{sato} for mild sufficient
conditions for this property to hold).
\end{itemize}

As it is usually done with L\'evy processes, we shall decompose $X$
into a compound Poisson process and a process with bounded jumps. More
specifically, for any $\varepsilon\in(0,1)$, we select a function
$c_{\varepsilon}\in C^{\infty}(\bbr)$, which is decreasing on $(-\infty
,0)$ and increasing on $(0,\infty)$ and such that $\mathbf{1}_{|x|\geq
\varepsilon}\leq c_{\varepsilon}(x)\leq\mathbf{1}_{|x|\geq \varepsilon
/2}$. Next, we define the truncated L\'evy densities
\[
{s_{\varepsilon}(x):= c_{\varepsilon}(x) s(x) \quad\mbox{and}\quad
\bar{s}_{\varepsilon}(x):= \bar{c}_{\varepsilon}(x) s(x),}
\]
with $\bar{c}_{\varepsilon}(x):=1- c_{\varepsilon}(x)$. Let
$Z^\varepsilon$ be a compound Poisson process with L\'evy measure
$s_\varepsilon(x) \,\mathrm{d}x$ and $X^\varepsilon$ be a L\'evy process, independent
from $Z^\varepsilon$, with characteristic triplet $(\sigma^2,\bar
s_\varepsilon(x) \,\mathrm{d}x,\mu_\varepsilon)$, where
%
%
\begin{equation}
\label{Dfnmuepsilon} {\mu_{\varepsilon}:=\mu-\int_{|x|\leq 1}x
c_{\varepsilon}(x)s(x)\,\mathrm{d}x}.
\end{equation}
It is clear that $X^\varepsilon+
Z^\varepsilon$ has the same law as $X$
and that the intensity and probability density of
{the} jumps of ${Z}^{\varepsilon}$ are $\lambda_{\varepsilon}:=\int
s_{\varepsilon}(x)\,\mathrm{d}x $ and
{$s_{\varepsilon}(x)/\lambda_{\varepsilon}$},
respectively. Throughout the paper, we let
$(N_{t}^{\varepsilon})_{t\geq 0}$ be the jump counting process of
$Z^{\varepsilon}$
and $(Y_{k}^{\varepsilon})_{k\geq 1}$ be the jump sizes of
$Z^{\varepsilon}$. Thus,
\(
Z_{t}^{\varepsilon}=\sum_{k=1}^{N_{t}^{\varepsilon}} Y_{k}^{\varepsilon}.
\)
Note that the distribution of $X_{t}^{\varepsilon}$ is also absolutely
continuous since $\sigma>0$ or $\int\bar{s}_{\varepsilon}(x)\,\mathrm{d}x=\infty
$, for any $\varepsilon>0$. For future reference, let us remark that
%
%
\begin{eqnarray}\label{FrmlMnVarSC}
\nonumber
\bbe \bigl(X^{\varepsilon}_{t} \bigr)&=&t \biggl(
\mu_{\varepsilon}+\int_{|x|\geq 1}x\bar{s}_{\varepsilon}(x)\,\mathrm{d}x
\biggr)=t\mu_{\varepsilon},
\nonumber
\\[-8pt]
\\[-8pt]
\nonumber
\operatorname{ Var} \bigl(X^{\varepsilon}_{t} \bigr)&=&t \biggl(
\sigma^{2}+\int x^{2}\bar{s}_{\varepsilon}(x)\,\mathrm{d}x \biggr)=: t
\sigma_{\varepsilon}^{2},
\end{eqnarray}
since $\varepsilon\in(0,1)$ (see, e.g., Example 25.12 in Sato \cite{sato}
for the mean and variance formulas of a L\'evy process).

The following lemma will be needed in what follows (cf. Propositions
I.4 and III.2 in L\'eandre~\cite{Leandre}). {See also Sections~\ref{SprmBnd}--\ref
{Subs:BndDstySmJmp} below for explicit expressions for the constants
$C_{p}(\eta,\varepsilon)$ and $c_{p}(\eta,\varepsilon)$.}
%
\begin{lmma}\label{KyLm1}
Let $f_{t}^{\varepsilon}$ be the transition density of the small-jump
component process $(X_{t}^{\varepsilon})_{t\geq 0}$. Then, for any
fixed positive real $\eta$ and positive integer $p$, there {exist} an
$\varepsilon_{0}(\eta,p)>0$ and {positive constants} $t_{0}(\eta
,\varepsilon)$, $c_{p}(\eta,\varepsilon)$, and $C_{p}(\eta,\varepsilon
)<\infty$ for any $\varepsilon<\varepsilon_{0}$ such that
%
%
\begin{eqnarray}
\label{KyInBn} \mathrm{(i)}\quad \bbp \Bigl(\sup_{0\leq s\leq t}\bigl|X_{s}^{\varepsilon}\bigr|
\geq \eta \Bigr)<C_{p}(\eta,\varepsilon) t^{p},\qquad \mathrm{(ii)}\quad
\sup_{|x|\geq \eta}f_{t}^{\varepsilon}(x)<c_{p}(
\eta ,\varepsilon) t^{p}
\end{eqnarray}
for all $0<t\leq t_{0}$ and $0<\varepsilon\leq \varepsilon_{0}$.
\end{lmma}

The following result provides the key tool for establishing the
small-time asymptotics of the moments of the L\'evy bridge ``stopped''
at the exit time from an interval $(a,b)$. Its proof is presented in
Appendix \ref{Sec:ProofTh:STBS}.
%
\begin{thrm}\label{Th:STBS}
For fixed constants $a\in[-\infty,0)$ and $b\in(0,\infty]$, define
\[
\tau:=\inf \bigl\{u\geq0\dvt X_{u}\notin(a,b) \bigr\}.
\]
Let $\varphi\dvtx \mathbb R\to\mathbb R$ be bounded and Lipschitz on $\bbr$
and let $\delta_{0}\in (0,\frac{b-a}{2} )$.
Then, for any $y\in(a+\delta_{0},b-\delta_{0})$ and {$0<\delta<\delta_{0}$},
%
%
\begin{eqnarray}
\label{Eq:FMR}
\bbe \bigl(\varphi(X_{\tau}){
\mathbf{1}}_{\{\tau\leq t, X_{t}\in({y-\delta
},y+\delta)\}} \bigr)=\int_{{y-\delta}}^{y+\delta} \biggl(
\frac{t^{2}}{2}\int_{(a,b)^{c}} \varphi(v)s(v)s(u-v)\,\mathrm{d}v+\calR
_{t}(u)t^{2} \biggr)\,\mathrm{d}u,\qquad
\end{eqnarray}
where the remainder term ${\calR_{t}(u)}$ is such that
%
%
\begin{equation}
\label{Eq:RmdBhv} \lim_{t\to 0}\mathop{\operatorname{ess\,sup}}_{u\in(a+\delta_{0},b-\delta_{0})}\bigl\llvert \mathcal{R}_{t}(u)\bigr\rrvert =0.
\end{equation}
\end{thrm}

\begin{remk}\label{EstExtProb}
By the definition of conditional expectation,
%
%
\begin{eqnarray}\label{DmbCER}
\nonumber
\bbe \bigl(\varphi(X_{\tau})\mathbf{1}_{\{\tau\leq t, X_{t}\in({y-\delta
},y+\delta)\}} \bigr)&=&
\bbe \bigl(\bbe \bigl(\varphi(X_{\tau})\mathbf{ 1}_{\{\tau\leq t\}}
\vert X_{t} \bigr)\mathbf{1}_{\{X_{t}\in({y-\delta
},y+\delta)\}} \bigr)
\nonumber
\\[-8pt]
\\[-8pt]
\nonumber
&=&\int_{{y-\delta}}^{y+\delta} \bbe \bigl(
\varphi(X_{\tau})\mathbf{1}_{\{\tau\leq t\}}\vert X_{t}=u \bigr)
f_{t}(u) \,\mathrm{d}u,
\end{eqnarray}
where $f_{t}(u)$ is the density of $X_{t}$ and, as usual,
$\Phi(u):=\bbe (\varphi(X_{\tau})\mathbf{
1}_{\{\tau\leq t\}}\vert  X_{t}=u )$ is such that
$\Phi(X_{t})$ is a version of $\bbe (\varphi(X_{\tau})\mathbf{
1}_{\{\tau\leq t\}}\vert  X_{t} )$. Comparing (\ref{DmbCER})
and (\ref{Eq:FMR}), it then follows that, for $\calL$-a.e. $y\in(a,b)$,
%
%
\begin{eqnarray}
\bbe \bigl(\varphi(X_{\tau})\mathbf{1}_{\{\tau\leq t\}}\vert
X_{t}=y \bigr)=\frac{{t^{2}}/{2}{\int_{(a,b)^{c}} \varphi
(v)s(v)s(y-v)\,\mathrm{d}v}}{f_{t}(y)}+\frac{\calR_{t}(y)t^{2}}{f_{t}(y)}.\label
{condovershoot}
\end{eqnarray}
If, in addition, the transition density $f_{t}$ satisfies the
asymptotic formula (\ref{Eq:AsymDstyPureLevy})\footnote{As stated in
the \hyperref[sec1]{Introduction}, (\ref{Eq:AsymDstyPureLevy}) holds for a large class of
Markov processes with jumps as proved by L\'eandre \cite{Leandre}. For
L\'evy processes, R\"uschendorf and
Woerner \cite{ruschendorf02} provided a more
elementary proof using the same conditions and similar
approach as in L\'eandre \cite{Leandre}. Higher order short-time
expansions for the transition densities were obtained in
Figueroa-L\'opez, Gong and
Houdr\'e \cite{FigGongHou:2011}.}
then, for $\calL$-a.e. $y\in(a,b)\setminus\{0\}$,
%
%
\begin{eqnarray}
\bbe \bigl(\varphi(X_{\tau})\mathbf{1}_{\{\tau\leq t\}}\vert
X_{t}=y \bigr) = t \frac{{\int_{(a,b)^{c}} \varphi(v)
s(v)s(y-v)\,\mathrm{d}v}}{2s(y)}+\mathrm{o}(t).\label{condovershootdens2}
\end{eqnarray}
Formulas \eqref{Eq:FMR} and \eqref{condovershoot}
can be interpreted as large deviation results for the trajectories of
L\'evy processes in small time.
When $\varphi(x)\equiv1$, (\ref{condovershootdens2}) gives the
following small-time approximation for the exit probability of the L\'
evy bridge:
%
%
\begin{equation}
\label{ExitProbaSmallTime} \bbp (\tau\leq t \vert X_{t}=y
)=t \frac{{\int_{(a,b)^{c}}s(v)s(y-v)\,\mathrm{d}v}}{2s(y)}+\mathrm{o}(t).
\end{equation}
\end{remk}

We conclude this section with a simpler result for the case when $X_t$
is outside the interval. Its proof is outlined in Appendix \ref
{Sec:ProofTh:STBS}.
%
\begin{prop}\label{Th:y_out}
Let $\varphi:\mathbb R\to\mathbb R$ be bounded and Lipschitz {on $\bbr
$}, and let $\delta_{0}>0$. Then, under the same notation and
conditions as in Theorem \ref{Th:STBS}, for any $y\in(a-\delta
_{0},b+\delta_{0})^c$ and $\delta<\delta_{0}$,
%
%
\begin{equation}
\label{Eq:FMRSmplCse} \bbe \bigl(\varphi(X_{\tau})\mathbf{1}_{\{X_{t}\in(y-\delta,y+\delta)\}
}
\bigr)=\int_{{y-\delta}}^{y+\delta} \bigl({t \varphi(u)s(u)}+
\calR_{t}(u)t \bigr)\,\mathrm{d}u,
\end{equation}
where the remainder term ${\calR_{t}(u)}$ is such that
%
%
\begin{equation}
\lim_{t\to 0}\mathop{\operatorname{ess\, sup}}_{u\in(a-\delta_{0},b+\delta_{0})^{c}}
\bigl\llvert \mathcal{R}_{t}(u)\bigr\rrvert =0.
\end{equation}
\end{prop}
%
\begin{remk}\label{EstExtProb2}
{Analogously to Remark \ref{EstExtProb}, (\ref{Eq:FMRSmplCse}) enables
us to establish the following natural asymptotic {formula}:
\[
\bbe \bigl(\varphi(X_{\tau})\vert X_{t}=y \bigr)=
t \frac{\varphi(y) s(y)}{f_{t}(y)}+{\mathrm{o}(1)}= \varphi(y) +\mathrm{o}(1)\qquad (t\to 0)
\]
for $\calL$-a.e. $y\in[a,b]^{c}$. The second equality above holds
whenever $f_{t}(y)$ satisfies (\ref{Eq:AsymDstyPureLevy}).}
\end{remk}

\section{On a precise bound for the remainder term}\label{estgen}
In the previous section, we developed the necessary results for finding
estimates of the functional
%
%
\begin{equation}
\label{Eq:FrstDfnKF} f(0,y,t):= \mathbb E \bigl[\varphi(X_\tau)
{\mathbf{1}}_{\tau\leq
t}\vert X_t = y \bigr]
\end{equation}
in short-time. Indeed, as explained in Remark \ref{EstExtProb},
Theorem \ref{Th:STBS} yields the following natural estimate for $f(0,y,t)$:
%
%
\begin{equation}
\label{Eq:FrstDfnEstKF} \tilde{f}(0,y,t)= \frac{{t^{2}}/{2}{\int_{(a,b)^c} \varphi(v)
s(v)s(y-v)\,\mathrm{d}v}}{f_{t}(y)}.
\end{equation}
The estimate {(\ref{Eq:FrstDfnEstKF})} will be used {below} to develop
adaptive discretization schemes for the Monte
Carlo computation of functionals of {the} killed L\'evy process {(see
Section~\ref{adaptive.sec})}. To this end, we first need to find an
explicit estimate for the remainder $\mathcal{R}_{t}(y)$ appearing in
(\ref{Eq:FMR}). Such an estimate can be expressed in terms of bounds
for the tail probability and transition densities of the small-jump
component $(X_{t}^{\varepsilon})_{t\geq 0}$. Hence, we start by
providing explicit expressions for the {upper bounds appearing} in (\ref
{KyInBn}) and then proceed to give {a} precise error bound for
$|f(0,y,t)-\tilde{f}(0,y,t)|$.

\subsection{Bounding the tail probability of the supremum}\label{SprmBnd}

The following exponential inequality for L\'evy processes with bounded
jumps {will be important to} estimate the supremum of the small-jump
component $(X_{t}^{\varepsilon})$ defined in Section~\ref{asymp}. Its proof, which is provided in Appendix
\ref{Ap:PrfErrBnd} for completeness, is a variation of {the} bound
{obtained} in R\"uschendorf and
Woerner \cite{ruschendorf02} (which in turn is based on Lemma
26.4 in Sato \cite{sato}).

\begin{lmma}\label{Lmm:KTPBnd}
Let $(M_t)$ be a martingale L\'evy process with $|\Delta M_t| \leq
\varepsilon$ and $\langle M,M\rangle_t = \sigma_\varepsilon^{2} t$. Then,
%
%
\begin{eqnarray}
\mathbb P \Bigl(\sup_{s\leq t} (M_s+ \mu s ) \geq
\eta \Bigr) \leq t^{{\eta}/{\varepsilon}}\bar{C}_{\ell}(\eta,\varepsilon;\mu)\qquad (
\ell={0,1}), \label{cgammaepsBar}
\end{eqnarray}
with the following constants {$\bar{C}_{\ell}(\eta,\varepsilon;\mu)$}
and corresponding conditions:
\begin{enumerate}[(1)]
\item[(1)] {$\bar{C}_{0}(\eta,\varepsilon;\mu)=
\mathrm{e}^{{\mu\vee0}/{\varepsilon}\mathrm{e}^{-1}}
\mathrm{e}^{\sigma_\varepsilon^2/\varepsilon^{2}}$ for all
$\eta>0$ and $0<t<\eta/(\mu\vee0)$} {(with the
convention here and below that the fraction is $+\infty$ if the
denominator is zero)};
\item[(2)] {$\bar{C}_{1}(\eta,\varepsilon;\mu)= \mathrm{e}^{{\mu\vee
0}/{\varepsilon}\mathrm{e}^{-1}} (\frac{\mathrm{e}\sigma_\varepsilon^2}{\varepsilon
\eta} )^{\eta/\varepsilon}$ for all $\eta>0$ and $0<t< \eta/(\mu
\vee0)$ if either \textup{(i)} $\mu\leq 0$ or \textup{(ii)} $\mu>0$} and $\eta\leq
\sigma_{\varepsilon}^{2}/\varepsilon$;
\end{enumerate}
\end{lmma}
In order to apply Lemma \ref{Lmm:KTPBnd} for $(X_{t}^{\varepsilon
})_{t\geq 0}$, we recall that $0<\varepsilon<1$ so that $\bbe
X_{t}^{\varepsilon}=\mu_{\varepsilon}t$. Then, the martingale part
{$M_{t}^{\varepsilon}:=X^{\varepsilon}-\mu_{\varepsilon}t$} of
$X^{\varepsilon}$ is such that
\[
{\bigl\anglel M^{\varepsilon},M^{\varepsilon}\bigr\angler _{t}=
\biggl(\sigma^{2}+\int \bar{c}_{\varepsilon}(x)x^{2}\nu(\mathrm{d}x)
\biggr) t=\sigma_{\varepsilon}^{2} t.}
\]
Thus, fixing
%
%
\begin{equation}
\label{Eq:DfnCnstEst} {t}_{0}(\varepsilon,\eta):= {
\frac{\eta}{2(\mu_{\varepsilon}\vee0)}},
\end{equation}
it follows that, for all $0<t<t_{0}$,
%
%
\begin{equation}
\label{Eq:AltMethBnd} \bbp \Bigl(\sup_{s\leq t}
X^{\varepsilon}_{s}\geq \eta \Bigr)\leq \bbp \Bigl(\sup
_{s\leq t} M^{\varepsilon}_{s}+|{\mu_{\varepsilon
}}|t
\geq \eta \Bigr)\leq \bbp \biggl(\sup_{s\leq t}
M^{\varepsilon}_{s}\geq \frac{\eta}{2} \biggr) \leq
t^{{{\eta}/{(2\varepsilon)}}}C \biggl(\frac{\eta}{2},\varepsilon \biggr),
\end{equation}
with $C(\eta,\varepsilon)$ is defined by
%
%
\begin{equation}
\label{cgammaeps} {C}(\eta,\varepsilon):= \biggl(\frac{\mathrm{e}\sigma_\varepsilon^2}{\varepsilon
\eta}
\biggr)^{{\eta}/{\varepsilon}}.
\end{equation}
Similarly, we have $\bbp(\sup_{s\leq t}
|X^{\varepsilon}_{s}|\geq \eta)\leq  2 t^{{{\eta}/{(2\varepsilon
)}}}C (\eta/2,\varepsilon )$.

\subsection{Bounding the transition density of the small-jump
component}\label{Subs:BndDstySmJmp}

{To obtain explicit expressions for the constants appearing in the
bounds for the density $f_{t}^{\varepsilon}$ in Lemma \ref{KyLm1},
we shall assume that the process $X$ is such that $X^\varepsilon_t$ has
a unimodal
distribution for all $t>0$ and $\varepsilon>0$.
By Yamazato's theorem (see Theorem 53.1 in Sato \cite{sato}), a sufficient
condition for this is that the process $X$ {is} self-decomposable,
which is the case if and only if the
L\'evy density $s$ is of the form\vspace*{1pt} $s(x) = \frac{k(x)}{|x|}$ for a
function $k$ which is increasing on $(-\infty,0)$ and decreasing on
$(0,\infty)$ (see Corollary 15.11 in Sato \cite{sato}). In particular, most
of the parametric models used in the literature (such as stable,
tempered stable, variance gamma, and normal inverse Gaussian processes)
are self-decomposable and so these processes as well as their
truncated versions have unimodal densities at all times. }

Let $m_t^{\varepsilon}$ be {the mode of $X_{t}^{\varepsilon}$}.
If $m_t^{\varepsilon} \in[-\underline\eta,\underline\eta]$ and $\eta
>\underline\eta$, then
the density can be estimated by
%
%
\begin{equation}
\label{Eq:Bnd1} \sup_{|x|\geq\eta} f^\varepsilon_t(x)
\leq\frac{2}{\eta-\underline
\eta} P\bigl[\bigl|X^\varepsilon_t\bigr|\geq\underline
\eta\bigr],
\end{equation}
simply because {the} density is decreasing in $(\underline\eta,\infty)$
and increasing in $(-\infty,-\underline\eta)$. The relation~(\ref
{Eq:Bnd1}) {in turn leads} to a bound of the form (\ref{KyInBn})(ii) by
applying the tail bound (\ref{KyInBn})(i).
It remains to find conditions for $m_t^{\varepsilon} \in[-\underline
\eta,\underline\eta]$. Since obviously $X^{\varepsilon}$ has finite
second moment, the following bound due to {Johnson and Rogers \cite
{JohnsonRogers}} can be applied
%
%
\begin{equation}
\label{JnsRgrsBnd} \bigl\llvert m_{t}^{\varepsilon}- \bbe
X_{t}^{\varepsilon}\bigr\rrvert ^{2}\leq3 \operatorname{ Var}
\bigl(X_{t}^{\varepsilon}\bigr).
\end{equation}
Thus, recalling the mean and variance formulas given in (\ref
{FrmlMnVarSC}), $m_{t}^{\varepsilon} \in[-\underline\eta,\underline\eta
]$ whenever $0<t<t_{1}$, where $t_{1}$ is such that
%
%
\begin{equation}
\label{Eq:Dfnt_1} t_{1} |{\mu_{\varepsilon}}|+\sqrt{3}
t_{1}^{1/2} \biggl(\sigma ^{2}+\int
_{|x|\leq\varepsilon} |x|^{2}\nu(\mathrm{d}x) \biggr)^{1/2}=
\underline \eta.
\end{equation}
By taking $\underline{\eta}=\eta/2$, we will have
%
%
\begin{equation}
\label{Eq:Bnd2} \sup_{|x|\geq\eta} f^\varepsilon_t(x)
\leq\frac{4}{\eta} \bbp \biggl[\bigl|X^\varepsilon_t\bigr|\geq
\frac{\eta}{2} \biggr]\leq \frac{8
C(\eta/4,\varepsilon)}{\eta} t^{{{\eta}/{(2\varepsilon)}}}
\end{equation}
for any {$0<t<t_{1}\wedge t_{0}$} with $t_{0}$ {defined} as in (\ref
{Eq:DfnCnstEst}).

\subsection{Precise bound for the remainder}

We are now ready to give an explicit bound for the reminder term
$\mathcal{R}_{t}(y)$ appearing in (\ref{Eq:FMR}), which in turn will
produce an error bound for $|f(0,y,t)-\tilde{f}(0,y,t)|$. Throughout,
we shall use the following notation:
\begin{longlist}[(iii)]
\item[(i)] $a_{\varepsilon}:=\sup_{x} s_{\varepsilon}(x)$ and
$a'_{\varepsilon}:=\sup_{x} |s'_{\varepsilon}(x)|$, where, as before,
$s_{\varepsilon}(x):=c_{\varepsilon}(x) s(x)$ is the L\'evy density
$s$, truncated in a neighborhood of the origin;
\item[{(ii)}] $\lambda_{\varepsilon}:=\int s(x)
c_{\varepsilon}(x)\,\mathrm{d}x$, ${
\mu_{\varepsilon}:=\mu}-\int_{|x|\leq 1}x
c_{\varepsilon}(x)s(x)\,\mathrm{d}x$, and
$\sigma_{\varepsilon}^{2}:=\sigma^{2}+\int
\bar{c}_{\varepsilon}(x) x^{2} {s(x)\,\mathrm{d}x}$;
\item[{(iii)}]
$C(\eta,\varepsilon)$ is defined as in
{\eqref{cgammaeps}}, $t_{0}(\varepsilon,\eta)$ is {defined} as in (\ref
{Eq:DfnCnstEst}), and $t_{1}(\varepsilon,\underline\eta)$ is {defined}
as in~(\ref{Eq:Dfnt_1}).
\end{longlist}

The following result, whose proof is given in Appendix \ref
{Ap:PrfErrBnd}, gives an estimate for $\mathcal{R}_{t}(y)$ in terms of
the previously defined notation and the $L_{\infty}$- and Lipschitz
norms of $\varphi$ denoted hereafter by
\begin{eqnarray*}
\|\varphi\|_{\infty}&:=&\operatorname{ ess}
\sup_{x}\bigl|\varphi(x)\bigr|,\\
\|\varphi\| _\mathrm{ Lip}&:=&\inf \bigl\{K\geq 0\dvt \bigl|\varphi(x)-\varphi(y)\bigr|\leq K
|x-y|, \forall x,y\in\bbr \bigr\}.
\end{eqnarray*}

\begin{thrm}\label{Th:KRFFEB}
Using the notation of Theorem \ref{Th:STBS}, assume that the process
$X$ is such that $X^\varepsilon_t$ has a unimodal
distribution for all $t>0$ and $\varepsilon>0$. Let $c:=b\wedge|a|$
and $\Delta_y:= (b-y)\wedge(y - a)>0$. Then,
\[
\bigl\llvert \mathcal{R}_{t}(y)\bigr\rrvert \leq{\frac{1}{t^{2}}}
e_{\calR}(0,y,t),\vadjust{\goodbreak}
\]
{for all $0<t<t_{0}(\varepsilon,(\Delta_y/2)\wedge c ) \wedge
t_{1}(\varepsilon,\Delta_y/2)$}, where
%
%
\begin{eqnarray}\label{TenOvBnd}
\nonumber
e_{\calR}(0,y,t)&:=&\mathrm{e}^{-\lambda_{\varepsilon}
t}\|\varphi\|_\infty
C(\Delta_y/4,\varepsilon) t^{{{\Delta
_y}/{(4\varepsilon)}}} \biggl\{
\frac{8}{\Delta_y} + 2a_\varepsilon t + a_\varepsilon
\lambda_\varepsilon t^2 \biggr\}
\\
&&{} + 2 \mathrm{e}^{-\lambda_{\varepsilon} t}\|\varphi\|_\infty a_{\varepsilon} C(c/2,
\varepsilon) t^{1+{{c}/{(2\varepsilon)}}} \{1+ {t \lambda_\varepsilon} \}
\nonumber
\\[-8pt]
\\[-8pt]
\nonumber
 &&{} +\frac{\|\varphi\|_\infty\lambda_{\varepsilon}^{2}
a_{\varepsilon}}{2} t^{3}+\|\varphi
\|_\infty a_{\varepsilon}\lambda _{\varepsilon}^{-1} \bigl(1
-\mathrm{e}^{-\lambda_{\varepsilon}t}\bigl[1+\lambda_{\varepsilon}t+ (\lambda_{\varepsilon}t)^{2}/2
\bigr] \bigr)
\\
&&{} + \mathrm{e}^{-\lambda_{\varepsilon}t}t^{2} \bigl[a_{\varepsilon}\lambda_{\varepsilon}
\|\varphi\|_\mathrm{ Lip}+{ 2}\|\varphi\|_{\infty}a_{\varepsilon}^{2}+
\|\varphi\|_{\infty}\lambda _{\varepsilon}a_{\varepsilon}'
\bigr]\biggl(\sigma_{\varepsilon}t^{1/2} + \frac{|{\mu_\varepsilon}|}{2} t\biggr).
\nonumber
\end{eqnarray}
\end{thrm}

Two immediate conclusions can be drawn. First, note that, by taking
$\varepsilon< {\frac{\Delta_y}{8} \wedge
\frac{c}{2}}$, we obtain a bound for the remainder satisfying condition
\eqref{Eq:RmdBhv}. Second, as seen in Remark \ref{EstExtProb}, the
previous {bound} implies the following error bound
\[
\bigl|f(0,y,t)-\tilde{f}(0,y,t)\bigr|\leq\frac{e_{\calR
}(0,y,t)}{f_{t}(y)}=:e_{f}(0,y,t),
\]
with $f$ and $\tilde{f}$ defined as in (\ref{Eq:FrstDfnKF})--(\ref
{Eq:FrstDfnEstKF}).

\begin{remk}\label{Rem:ApprxExitProba}
The approximation for the conditional exit probability
{$p(0,y,t):= \mathbb P [\tau\leq t\vert X_t = y ]$} is
obtained by substituting $\varphi\equiv1$ into \eqref{condovershoot}:
\[
\tilde{p}(0,y,t)= \frac{{t^{2}}/{2}{\int_{(a,b)^c}
s(v)s(y-v)\,\mathrm{d}v}}{f_{t}(y)}.
\]
Making this substitution in the previous bound, it follows that
$|p(0,y,t)-\tilde{p}(0,y,t)|\leq e_{p}(0,y,t)$ with $e_{p}(0,y,t)$
given by
%
\begin{eqnarray*}
e_{p}(0,y,t)&:=&\frac{1}{f_{t}(y)} \biggl(\mathrm{e}^{-\lambda_{\varepsilon}
t} C(
\Delta_y/4,\varepsilon) t^{{{\Delta_y}/{(4\varepsilon)}}} \biggl\{\frac{8}{\Delta_y}
+ 2a_\varepsilon t + a_\varepsilon \lambda_\varepsilon
t^2 \biggr\}
\nonumber\\
&&\hspace*{28pt}{} + 2\mathrm{e}^{-\lambda_{\varepsilon} t}a_{\varepsilon} C(c/2,\varepsilon) t^{1+{{c}/{(2\varepsilon)}}} \{1+
{t \lambda_\varepsilon} \} +\frac{\lambda_{\varepsilon}^{2}
a_{\varepsilon}}{2} t^{3}
\nonumber
\\[-8pt]
\\[-8pt]
\nonumber
&&\hspace*{28pt}{} + a_{\varepsilon}\lambda_{\varepsilon}^{-1} \bigl(1
-\mathrm{e}^{-\lambda_{\varepsilon}t}\bigl[1+\lambda_{\varepsilon}t+ (\lambda_{\varepsilon}t)^{2}/2
\bigr] \bigr)
\\
&&\hspace*{28pt}{} + \mathrm{e}^{-\lambda_{\varepsilon}t}t^{2} \bigl[{ 2}a_{\varepsilon}^{2}+
\lambda_{\varepsilon}{a'_{\varepsilon}} \bigr]\biggl(
\sigma_{\varepsilon}t^{1/2} + \frac{|{\mu_\varepsilon}|}{2} t\biggr) \biggr),
\nonumber
\end{eqnarray*}
valid for all $t<t_{0}(\varepsilon,(\Delta_y/2)\wedge c ) \wedge
t_{1}(\varepsilon,\Delta_y/2)$. The one-sided case ($a=-\infty$) can
similarly be obtained.
\end{remk}

\section{Adaptive simulation of killed L\'evy processes}
\label{adaptive.sec}

Our goal in this section is to design {a type of} adaptive Monte Carlo
estimators for {functionals} of the form
%
%
\begin{equation}
\mathbb E\bigl[F(X_T){\mathbf{1}_{\tau> T}}\bigr],
\label{killed}
\end{equation}
where ${F}$ is a Borel measurable function and $\tau:= \inf\{t\geq0\dvt X_t \notin D\}$ with $D := (a,b)$, for some $a\in[-\infty,0)$ and
$b\in(0,\infty]$.
From now on, to simplify notation and with no loss of generality, we
shall take $T=1$.

For {$0<s<t$}, $x\in\mathbb R$, and $y\in\mathbb R$, we denote by
$\mathbb P^{\mathrm{BR}}_{(s,t,x,y)}[\cdot]$
the bridge law of the L\'evy process~$X$ on the time interval $[s,t]$ with
starting value $x$ and terminal value $y$; that is, this is a version
of the regular
conditional distribution of 
{$\{x+X_{u-s}\}_{u\in[s,t]}$ given $X_{t-s}=y-x$}.
Since $X_t$ has a strictly
positive density on $\mathbb R$ for every $t>0$, the bridge law is
uniquely defined for
$\calL$-almost every $y\in\mathbb R$ (recall that $\calL$ stands for
the Lebesgue measure), which is sufficient for our purposes. We also
let $p(x,y,t)$ denote the exit probability from the domain $D$ before
time $t$ for the L\'evy bridge:
%
%
\begin{equation}
\label{Eq:Dfnp} p(x,y,t):= \mathbb P^{\mathrm{BR}}_{(0,t,x,y)}[\tau\leq t] =
\bbp \bigl[\exists u\in[0,t]\dvt x+X_{u}\notin(a,b)\vert
X_{t}=y \bigr].
\end{equation}

Our approach is based on the following decomposition:
%
%
\begin{equation}
\label{Eq:DtrmTmCs} \mathbb E\bigl[{F} (X_1){\mathbf{1}_{\tau> 1}}
\bigr]=\bbe \Biggl[{F} (X_1)\prod_{i=0}^{N-1}
\bigl(1-{p}(X_{T_{i}},X_{T_{i+1}},T_{i+1}-T_{i})
\bigr) \Biggr],
\end{equation}
where $0=T_{0}\leq\cdots\leq T_{N}=1$ are suitable sampling
times. Formula (\ref{Eq:DtrmTmCs}) {directly follows from the Markov
property} when the sampling points are deterministic. In that case,
the set of points $\calX:=\{(T_{i},X_{T_{i}})\}_{i=0}^{N}$ is called a
deterministic skeleton. In our setting, both the number of
points~$N$ and the sampling times $0=T_{0}\leq T_{1}\leq\cdots\leq
T_{N}=1$ are random and we need to formalize under what conditions on
$\calX$ (\ref{Eq:DtrmTmCs}) still holds. The following result will
suffice for our purposes.
%
\begin{lmma}\label{Lm:Skltn}
Let $N$ be a random variable with support
$\mathcal N \subseteq\mathbb N$, such that $N>0$, and let
$0=T_{0}\leq\cdots\leq T_{N}=1$ be random points such that
\begin{enumerate}[(1)]
\item[(1)] Each $T_{i}$ takes values in a {countable} set $\calK\subset[0,1]$;
\item[(2)] For each $n\in\calN$ and
$(s_{0},\ldots,s_{n})\in{\calK^{n+1}}$ with $0=s_{0}\leq \cdots\leq
s_{n}=1$, the event $\{N=n, (T_{0},\ldots,T_{n})=(s_{0},\ldots,s_{n})\}$
is $\sigma(X_{s_{i}}\dvt i=0,\ldots,n)$-measurable.
\end{enumerate}
Then, (\ref{Eq:DtrmTmCs}) is satisfied for any
measurable function $F$ with $\bbe[|{F}(X_1)|]<\infty$ {and,
furthermore,} for every $t\in(0,1)$, $n\in\mathcal N$, and $A\in
\mathcal B(\mathbb R)$,
%
%
\begin{eqnarray}
\label{condbridgelaw} \mathbb P [X_t \in A| N=n, T_0,
\ldots,T_N, X_{T_0},\ldots ,X_{T_N} ] = \mathbb
P^{\mathrm{BR}}_{T_{i^*}, T_{i^*+1}, X_{T_{i^*}},X_{T_{i^*}+1}} [X_t \in A ],
\end{eqnarray}
where $i^* = \max\{i\dvt T_i \leq t\}$.
\end{lmma}
\begin{pf}
Throughout, we let $\bar{p}(x,y,t):=1-p(x,y,t)$, $\vec\calK^{n}:=\{
(s_{0},\ldots,s_{n})\in{\calK^{n+1}}\dvt 0=s_{0}\leq \cdots\leq s_{n}=1\}$,
{$U_{0}:=[s_{0},s_{1}]$}, and $U_{i}:=(s_{i},s_{i+1}]$, with $i=1,\ldots
,n-1$. We also use the notation
%
%
\begin{equation}
\label{Eq:NdIndFc} \calI_{U}:=\mathbf{1}_{\{X_{u}\in(a,b)\dvt \forall  u\in U\}}\qquad\mbox{for a
domain }U\subset\bbr_{+} \mbox{ and } \calI_{\varnothing}=1.
\end{equation}
Then, by Markov property
\begin{eqnarray*}
&&\mathbb E\bigl[{F} (X_1){\mathbf{1}_{\tau> 1}}\bigr]\\
 &&\quad=\sum
_{n\in\calN}\sum_{(s_{0},\ldots,s_{n})\in\vec\calK^{n}} \mathbb
E \bigl[{F} (X_1) \calI_{[0,1]}\mathbf{1}_{\{{N}=n,(T_{0},\ldots
,T_{n})=(s_{0},\ldots,s_{n})\}} \bigr]
\\
&&\quad=\sum_{n\in\calN}\sum_{(s_{1},\ldots,s_{n})\in\vec\calK^{n}}
\mathbb E \Biggl[{F} (X_1)\mathbf{1}_{\{{N}=n,(T_{0},\ldots,T_{n})=(s_{0},\ldots
,s_{n})\}}\bbe \Biggl[\prod_{i=0}^{n-1}
\calI_{U_{i}} |X_{s_{j}}\dvt j=0,\ldots,n \Biggr] \Biggr]
\\
&&\quad=\sum_{n\in\calN}\sum_{(s_{1},\ldots,s_{n})\in\vec\calK^{n}}
\mathbb E \Biggl[{F} (X_1)\mathbf{1}_{\{N=n,(T_{0},\ldots,T_{n})=(s_{0},\ldots
,s_{n})\}}\prod_{i=0}^{n-1} \bar {p}
(X_{T_{i}},X_{T_{i+1}},T_{i+1}-T_{i} )
\Biggr]
\\
&&\quad =\mathbb E \Biggl[{F}(X_1)\prod_{i=0}^{N-1}
\bar{p} (X_{T_{i}},X_{T_{i+1}},T_{i+1}-T_{i}
) \Biggr],
\end{eqnarray*}
which proves (\ref{Eq:DtrmTmCs}). Similarly, $\mathbb P[X_t \in A| N=n,
T_0,\ldots,T_N, X_{T_0},\ldots,X_{T_N}]$ can be decomposed as
\begin{eqnarray*}
&&\sum_{(s_{0},\ldots,s_{n})\in\vec\calK^{n}}\mathbb P \bigl[X_t \in A|
N=n, (T_0,\ldots,T_n) = (s_0,
\ldots,s_n), X_{T_0},\ldots,X_{T_N} \bigr]
\mathbf1_{ (T_0,\ldots,T_n) = (s_0,\ldots,s_n)}
\\
&&\qquad = \sum_{(s_{0},\ldots,s_{n})\in\vec\calK^{n}} \mathbb P^{\mathrm{BR}}_{s_{i^*}, s_{i^*+1}, X_{s_{i^*}},X_{s_{i^*}+1}}[X_t
\in A]\mathbf1_{(T_0,\ldots,T_n) = (s_0,\ldots,s_n)}
\\
&&\qquad = \mathbb P^{\mathrm{BR}}_{T_{i^*}, T_{i^*+1},
X_{T_{i^*}},X_{T_{i^*}+1}}[X_t \in A].
\end{eqnarray*}
\upqed\end{pf}
From (\ref{Eq:DtrmTmCs}), it is now evident that, for the computation
of \eqref{killed} by Monte Carlo, it suffices to simulate independent
replicas of the random variable $\mathcal{Y}:=F(X_{1})N(\calX)$, where
hereafter we denote
\[
N(\mathcal X):= \prod_{i=0}^{N-1}
\bigl(1-p(X_{T_{i}},X_{T_{i+1}},T_{i+1} -
T_{i}) \bigr).
\]
The exit probability $p(x,y,t)$ does not typically admit
{a} closed form expression and some type of approximation must be
applied for
its evaluation. The short-time asymptotics
(\ref{condovershoot}) yields the following natural estimate for
$p(x,y,t)$ {when $x,y\in D$:
%
%
\begin{equation}
\label{Eq:DfnApproxExitProba1} \tilde{p}(x,y,t):=\bigl(\breve{p}(x,y,t)\vee0
\bigr)\wedge1\qquad \mbox {with } \breve{p}(x,y,t):= 
\frac{t^{2}}{2} \int_{(a,b)^{c}}\frac{s(u-x)s(y-u)}{f_{t}(y-x)} \,\mathrm{d}u.
\end{equation}
{We also set} $\tilde{p}(x,y,t)=1$ if $x\notin D$ or
$y\notin D$.} This approximation satisfies
%
%
\begin{eqnarray}
\bigl|\tilde p(x,y,t) - p(x,y,t) \bigr|\leq e_p(x,y,t), \label{exitest}
\end{eqnarray}
where the error bound $e_p(x,y,t)$ is defined as in Remark
\ref{Rem:ApprxExitProba} for $x,y\in D$ and by $e_p(x,y,t)=0$ if
$x\notin D$ or $y\notin D$. We can then approximate {$N(\mathcal X)$} by
%
%
\begin{equation}
\tilde N(\mathcal X):= \prod_{i=0}^{N-1}
\bigl(1-\tilde p(X_{T_{i}},X_{T_{i+1}},T_{i+1} -
T_{i}) \bigr). \label{Eq:ApprxProdMC}
\end{equation}

Replacing the true exit probability $p(x,y,t)$ with its approximation
$\tilde p(x,y,t)$ introduces a bias into the evaluation of $N(\mathcal
X)$, which is {hard} to quantify if the process
$X$ is discretized using the uniformly
spaced grid $T_{i}=i/N$.
For this reason, we now propose an adaptive algorithm for the
determination of the sampling times, which starts by simulating
the terminal value $X_1$ and then refines the sampling grid, using more
discretization points when the estimate of the approximation error is ``large''.
The algorithm is parameterized by a real number ${\gamma} >0$, which
represents the error tolerance and ensures
that under suitable conditions on $e_p$, the global discretization
error for approximating the quantity of interest
\eqref{killed} will be bounded by ${\gamma}$
(see Proposition \ref{biasprop} below). The algorithm also requires
simulation from the marginal distribution $f_{1}$ of $X_{1}$ and the bridge
distribution of $X_{t/2}$ conditioned to
$X_{t}=y$ ($t>0$). Hereafter, we denote the density of this bridge
distribution by ${f}^{\mathrm{br}}_{t/2}(x,y)$ and recall the following
{well-known} formula:
%
%
\begin{equation}
\label{Eq:BrdDsity1} {f}^{\mathrm{br}}_{t/2}(x,y) :=
\frac{{f}_{t/2}(x){f}_{t/2}(y-x)}{{f}_t(y)}.
\end{equation}
At the end of this section, we introduce a new method to simulate
variates from the density (\ref{Eq:BrdDsity1}).

\begin{algorithm}[t]
\caption{$[\calX] = \mathrm{GenerateSkeleton}(\gamma)$}
\floatname{algorithm}{Procedure}
\label{ConstrSkeleton}
\begin{algorithmic}
\STATE$N_{0}=0$, $N_{1}=1$, $m=1$
\STATE$T_{0}^{1}=0$, $T_{1}^{1}=1$, $X_0 = 0$
\STATE Generate an observation $X_1$ from the density $f_{1}$
\WHILE{$N_{m}\neq N_{m-1}$ {and $\{X_{T^m_i} \in D,\mbox{ for
}i=1,\ldots,N_m\}$}}
\STATE$n=0$, $T_{0}^{m+1}=0$
\FOR{$i=0 \to N_{m}-1$}
\STATE$\Delta T=T_{i+1}^{m}-T_{i}^{m}$
\IF{$e_{p}(X_{T_{i}^{m}},X_{T^{m}_{i+1}},\Delta T)> \gamma\Delta T$}
\STATE$T^{m+1}_{n+1}=(T_{i}^{m}+T_{i+1}^{m})/2$,
$T^{m+1}_{n+2}:=T_{i+1}^{m}$
\STATE Generate an observation $X_{T^{m+1}_{n+1}}$ from the bridge
density $f_{\Delta T/2}^{\mathrm{br}}(\cdot,X_{T_{i+1}^{m}}-X_{T_{i}^{m}})$
\STATE$n=n+2$
\ELSE
\STATE$T_{n+1}^{m+1}:=T_{i+1}^{m}$
\STATE$n=n+1$
\ENDIF
\ENDFOR
\STATE{$N_{m+1}=n$}
\STATE{$m=m+1$}
\ENDWHILE
\STATE RETURN $\calX=\{(T_{i}^{m},X_{T_{i}^{m}})\}_{i=0}^{N_{m}}$.
\end{algorithmic}

\end{algorithm}

The procedure to generate the skeleton of $X$ is outlined in
pseudo-code in Algorithm \ref{ConstrSkeleton}
below. Assume that this algorithm terminates in
finite time a.s. (see Proposition \ref{biasprop} for sufficient
conditions for this to hold). The algorithm then defines a pair $N$ and
$\calT:=(T_{0},\ldots,T_{N})$, which satisfies the conditions of Lemma
\ref{Lm:Skltn}. Indeed,
by construction,
each $T_{i}$ takes values in the dyadic grid
$\{i2^{-m},i=0,\ldots,2^{m},m=0,1,\ldots\}$, which is a countable set.
To check the second condition of the lemma, we fix $n$ and a partition
$\pi:=\{s_0,\ldots,s_n\}$ of $[0,1]$, and proceed as follows to write
the event $E:=\{N=n,T_{0}=s_{0},T_{1}=s_{1},\ldots,
T_{n}=s_{n}\}$ in terms of $\{X_{s_{i}}\}_{i=0}^{n}$:
\begin{itemize}
\item We can and will assume with no loss of generality that
$\pi$ is a \emph{recursive dyadic partition}, meaning that {$\{0,1\}
\subset
\pi$ and, for every $t\in(0,1)\cap\pi$,} there exists $k\in\mathbb
N$ with $2^k t \in\mathbb N$, and if we
take the smallest such $k$ then also $t+ \frac{1}{2^k} \in\pi$ and
$t- \frac{1}{2^k} \in\pi$. By construction, if $\pi$ does not have
this property, the event $E$ has zero probability.
%
\item We shall assume that $n\geq2$ because if $n=1$ then necessarily
$s_0 = 0$ and $s_1=1$ and, therefore,
\[
E = {\{X_1 \notin D\}\cup\bigl\{X_1 \in D,
e_p(X_0,X_1,1)\leq\gamma\bigr\}} \in
\sigma(X_0,X_1).
\]
\item For each $\ell\in\{0,\ldots,n-1\}$, define
$\pi_{\ell}:=\{s_{i}\in\pi\dvt  2^{n-\ell} s_{i} \mbox{ is an even
integer}\}$. The number of elements of $\pi_{\ell}$ is denoted
$n_{\ell}$ and the sorted elements of $\pi_{\ell}$ are denoted
$s_{1}^{\ell}<\cdots<s_{n_{\ell}}^{\ell}$. Clearly, $\pi_0 = \pi$ and
$\pi_{n-1} \neq\pi$ since $1/2\in\pi$ whenever $n\geq2$; we let $\ell
^* = \max\{l\geq0\dvt \pi_l = \pi\}$ and
$\pi^* = \pi\setminus\pi_{\ell^*+1}$.

\item For each $i=1,\ldots,n_{\ell}-1$, define the event
\[
E^{\ell}_{i}:= \bigl\{\omega\dvt  e_{p}
\bigl(X_{s_{i}^{\ell}}(\omega),X_{s_{i+1}^{\ell
}}(\omega),s_{i+1}^{\ell}-s_{i}^{\ell}
\bigr)\leq\gamma\bigl(s_{i+1}^{\ell
}-s_{i}^{\ell}
\bigr) \bigr\}
\]
if $\pi\cap(s_{i}^{\ell},s_{i+1}^{\ell})=\varnothing$; otherwise, we set
\[
E^{\ell}_{i}:= \bigl\{\omega\dvt e_{p}
\bigl(X_{s_{i}^{\ell}}(\omega ),X_{s_{i+1}^{\ell}}(\omega),s_{i+1}^{\ell}-s_{i}^{\ell}
\bigr)>\gamma \bigl(s_{i+1}^{\ell}-s_{i}^{\ell}
\bigr) \bigr\}.
\]
Then it follows that
\begin{eqnarray*}
E&=& \Biggl\{\bigcap_{i=0}^n
\{X_{s_i} \in D\} \cap\bigcap_{\ell=\ell
^*}^{n-1}
\bigcap_{i=1}^{n_{\ell}-1} E_{i}^{\ell}
\Biggr\}
\\
&&{} \cup \Biggl\{\bigcup_{s\in\pi^*} \{X_{s}
\notin D\} \cap\bigcap_{s\in\pi_{\ell^*+1}} \{X_{s} \in
D\} \cap\bigcap_{\ell=\ell^*+1}^{n-1}\bigcap
_{i=1}^{n_{\ell}-1} E_{i}^{\ell} \Biggr
\},
\end{eqnarray*}
which clearly belongs to $\sigma(X_{s_{i}}\dvt i=0,\ldots,n)$.
\end{itemize}
To see that $X_{T^{m+1}_{n+1}}$ can be sampled
from the bridge density $f^{\mathrm{br}}_{\Delta T/2} (\cdot, X_{T^m_{i+1}} -
X_{T^m_i})$ {in Algorithm~\ref{ConstrSkeleton}}, we can apply the
second part of Lemma
\ref{Lm:Skltn} to the couple {$(k, \mathcal T_k)$}, where
$\mathcal T_k = \{T_0,\ldots,T_k\}$ contains the first {$k+1$ sampling times}
which have been added to the grid by the algorithm, in increasing order.

Algorithm \ref{ConstrSkeleton} terminates when {at least one of the
sampling observations $X_{T_{i}}$ is out of the domain $D$ or} the
error over each subinterval of the sampling times $0=T_{0}<\cdots
<T_{N}=1$ is small enough in the following sense:
%
%
\begin{eqnarray}
e_p(X_{T_{i}},X_{T_{i+1}},T_{i+1}-T_{i})
\leq{\gamma} (T_{i+1}-T_{i}),\qquad i=0,\ldots,N-1.
\label{errcontrol}
\end{eqnarray}
{At} first glance, it is not obvious that the algorithm will actually
terminate in finite time.
The following result gives conditions under which {this is the case}
and shows that the global error of the estimate is of order $\gamma$.
%
\begin{prop}\label{biasprop}
The following assertions hold:
\begin{enumerate}[(ii)]
\item[(i)] Let $X$ be a L\'evy process
satisfying one of the {following} two (non-mutually exclusive) conditions:
\begin{enumerate}[1.]
\item[1.]$X$ does not hit points{;} that is, {$\bbp(\tau^{\{x\}}<\infty
)=0$ for all $x$, where $\tau^{\{x\}}:=\inf\{s>0\dvt X_{s}=x\}$ or, equivalently},
\[
\int_{\mathbb R}\Re \biggl(\frac{1}{1+\psi(u)} \biggr) \,\mathrm{d}u = \infty,
\]
where $\psi(u) = \log\mathbb E[\mathrm{e}^{\mathrm{i}uX_1}]$ {(see Kyprianou \cite{kyprianou}, Theorem
7.12)};
\item[2.]$X$ is a finite variation process.
\end{enumerate}
{Additionally,} assume that the upper bound of the approximation error
$e_p(x,y,t)$ satisfies
%
%
\begin{equation}
\lim_{t\downarrow0} \frac{1}{t} \sup_{x,y \in(a',b')}
e_p(x,y,t) = 0\qquad\forall a',b' \in(a,b).
\label{errbound}
\end{equation}
Then, Algorithm \ref{ConstrSkeleton} terminates in finite time a.s.
\item[(ii)] Assume that {$\bbe|F(X_1)| < \infty$}. Let $\mathcal X=\{
(T_{i},X_{T_{i}})\}_{i=0}^{N}$ be a
skeleton of $X$ on $[0,1]$ satisfying \eqref{errcontrol} and $\tilde
\calN(\calX)$ be given by (\ref{Eq:ApprxProdMC}). Then,
%
%
\begin{eqnarray}
\bigl|\bbe\bigl[{F}(X_1){\mathbf{1}}_{\tau>1}\bigr] - \bbe
\bigl[{F}(X_1) \tilde N(\mathcal X)\bigr]\bigr| \leq{\gamma} \bbe
\bigl[\bigl|{F}(X_1)\bigr|\bigr].\label{algoerr}
\end{eqnarray}
\end{enumerate}
\end{prop}

\begin{remk}\label{betterbias}
In view of {Proposition \ref{biasprop}}, $\bbe[{F}(X_1){\mathbf{1}}_{\tau>1}]$
can be approximated {by} the Monte Carlo estimator
\[
\frac{1}{M} \sum_{{k=1}}^M {F}
\bigl(X_1^{({k})} \bigr) \tilde N \bigl(\mathcal
X^{({k})} \bigr),
\]
where $X^{({k})}$ are {independent} copies of the process $X$ and
$\tilde N
(\mathcal X^{({k})})$ are corresponding values {computed with formula
(\ref{Eq:ApprxProdMC})}. This estimator has a statistical error which
can be
estimated in the usual way, and a discretization bias, which is
bounded from above by ${\gamma} \bbe[|{F}(X_1)|]$.
{In view of \eqref{precisebias} {below}, a more precise a posteriori
estimate of the bias
is
\[
\frac{1}{M}\sum_{k=1}^M\bigl|{F}
\bigl(X^{(k)}_1\bigr)\bigr| {\mathbf{1}}_{S^{(k)}_N} \sum
_{i=1}^{N} e_p \bigl(X^{(k)}_{T^{(k)}_{i}},X^{(k)}_{T^{(k)}_{i+1}},
T^{(k)}_{i+1}-T^{(k)}_{i} \bigr),
\]
with $S_N:=\{(X_{T_0},\ldots,X_{T_N})\in D^{N+1}\}$.}
\end{remk}

\begin{lmma}\label{jumplemma}
Let $X$ be a L\'evy process such that for all $t>0$, the law of $X_t$
has no atom. Then, for all $x\in\mathbb R$,
\begin{eqnarray*}
\mathbb P\bigl[\bigl\{t\in[0,1]\dvt \Delta X_t \neq0,
X_{t-} = x\bigr\} =\varnothing\bigr] = 1; \qquad\mathbb P\bigl[\bigl\{t
\in[0,1]\dvt \Delta X_t \neq0, X_{t} = x\bigr\} =\varnothing
\bigr] = 1.
\end{eqnarray*}
\end{lmma}
\begin{pf}
We only prove the first identity, the second one follows by similar
arguments (or alternatively by time reversal).
Let $N^\varepsilon_1 = \#\{t\in[0,1]\dvt |\Delta X_t| > \varepsilon,
X_{t-} =
x\}$.
Then
\[
\mathbb P\bigl[\bigl\{t\in[0,1]\dvt \Delta X_t \neq0,
X_{t-} = x\bigr\} \neq\varnothing\bigr] \leq\mathbb E
\bigl[N^0_1\bigr] \leq\sum_{n=1}^\infty
\mathbb E\bigl[N^{{1}/{n}}_1\bigr].
\]
But by the compensation formula (see Bertoin \cite{bertoin}, Section~0.5),
\[
\mathbb E\bigl[N^{{1}/{n}}_1\bigr] = \mathbb E \biggl[\int
_0^1 \int_{|y|>\varepsilon}
\mathbf{1}_{X_s = x} \nu(\mathrm{d}y) \,\mathrm{d}s \biggr] = \int_{|y|>\varepsilon} \nu(\mathrm{d}y)
\int_0^1 \mathbb P[X_s = x] \,\mathrm{d}s =
0.
\]
\upqed\end{pf}

\begin{pf*}{Proof of Proposition \ref{biasprop}}
\emph{Part} (i).  With the aim of obtaining a
contradiction, assume that the statement of the proposition is not true,
and the algorithm does not terminate. Let $\{\tilde T_i\}_{i\geq1}$
{be the infinite {sequence} of different sampling times} produced by
the algorithm (in the order in which they were generated, that is,
not necessarily ordered in time). {Let $\tilde X_i:=X_{\tilde{T}_{i}}$
be the corresponding sampling observations.}
Since the sequence $\{\tilde T_i\}$ is bounded, we
can find indices $\{{i_k}\}_{k\geq1}$ such that $\tilde T_{i_k} \to
T^*$. Moreover, since every point $\tilde T_i$ (for $i\geq2$) is
obtained as a
midpoint of a certain interval, we can find two sequences
$\{T^-_i\}$ and $\{T^+_i\}$ such that $T^-_i \uparrow T^*$, $T^+_i
\downarrow T^*$, $T^*\in[T^-_i,T^+_i]$ for all $i$ and
$e_p(X_{T_i^-},X_{T_i^+},T_i^+-T_i^-)> {\gamma}(T_i^+-T_i^-)$ for
all $i$. In addition, since the process $X$ has right and left
limits, both $\lim X_{T^+_i} = X^+$ and $\lim X_{T^-_i} = X^-$ {exist}.
There are three possibilities.

If $X^-\in(-\infty,a)\cup(b,\infty)$ or $X^+\in(-\infty,a)\cup
(b,\infty)$ then for some $i$,
$\tilde X_{i} \notin D$, so that the algorithm must have stopped
in finite time and we have a contradiction.

If $X^- \in(a,b)$ and $X^+ \in(a,b)$ then, using the
property {\eqref{errbound}}, we can find a contradiction with
$e_p(X_{T_i^-},X_{T_i^+},T_i^+-T_i^-)> {\gamma}(T_i^+-T_i^-)$.

It remains to treat the case when $X^-$ or $X^+$, or both, are at
the boundary of $D$.
Then, either $X^- = X^+ = X_{T^*}$ or $\Delta X_{T^*}\neq0$. The
latter case is ruled out by Lemma \ref{jumplemma} and in the case when
$X$ cannot hit points, the former case is ruled out as well.

We may
therefore assume that $X$ is a finite variation process with nonzero
drift $\mu$ ({cf.} Kyprianou \cite{kyprianou}, Theorem 7.12) and, to fix the notation,
that $X^- = X^+ = X_{T^*}=b$. We may also assume that $T^*$ is irrational,
since for every $t\in\mathbb Q \cap[0,1]$, {$\mathbb P[X_t = b] =
0$}. The fact that $T^*\notin\mathbb Q$ implies that $T_i^- <
T^* < T_i^+$ for every $i$, and we can also assume that $X_{T^+_i}$
and $X_{T^-_i}$ belong to $D$ for each $i$, because otherwise the
algorithm would have stopped in finite time.

Introduce two sequences
of stopping times:
\[
\sigma_{n} := \inf\{t>\tau_n\dvt X_t \leq b\}
\wedge1, \qquad \tau_{n+1} := \inf\{t>\sigma_{n}\dvt X_t
\geq b\}\wedge1, n\geq0,
\]
with $\tau_0 := \inf\{t>0\dvt X_t \geq b\}\wedge1$.
The sequences $\{\tau_n\}$ and
$\{\sigma_n\}$ do not have an accumulation point except $t=1$ and
for each $n\geq0$, $\sigma_{n} > \tau_n$ if $\tau_n<1$ and
$\tau_{n+1} > \sigma_n$ if $\sigma_n<1$. This holds because for a finite
variation process $X$ with drift $\mu\neq0$, $\{0\}$ is irregular
for $[0,\infty)$ if $\mu<0$ and for $(-\infty,0]$ if $\mu>0$
(Sato \cite{sato}, Theorem 43.20), and
$X$ may only creep in the direction opposite to the drift
(Kyprianou \cite{kyprianou}, Theorem 7.11). Then
clearly, for every $\tau\in[0,1]$ such that
$X_\tau= b$, either there is $n\geq0$ with $\sigma_n=\tau$, which
means that for some $\varepsilon>0$, $X_t \notin D$ for $t\in
(\tau-\varepsilon,\tau)$, or there is $n\geq0$ with
$\tau_n=\tau$, which
means that for some $\varepsilon>0$, $X_t \notin D$ for $t\in
(\tau,\tau+\varepsilon)$. In both cases, there is a contradiction
with the fact that $X_{T^+_i}$
and $X_{T^-_i}$ belong to $D$ for each $i$.

\emph{Part} (ii).  Below, we denote $\bar{p}(x,y,t):= 1-p(x,y,t)$,
$\bar{\tilde{p}}(x,y,t)=1-\tilde{p}(x,y,t)$, {and $S_N:=\{(X_{T_0},\ldots
,X_{T_N})\in
D^{N+1}\}$. Then, since}
\[
N(\calX)-\tilde{N}(\calX)=\prod_{i=0}^{N-1}
\bar{p} (X_{T_{i}},X_{T_{i+1}},T_{i+1}-T_{i}
)-\prod_{i=0}^{N-1} \bar {\tilde{p}}
(X_{T_{i}},X_{T_{i+1}},T_{i+1}-T_{i} ),
\]
we get
%
%
\begin{eqnarray}\label{precisebias}\hspace*{-30pt}
\bigl|\bbe\bigl[{F}(X_1){\mathbf{1}}_{\tau>1}\bigr] - \bbe
\bigl[{F}(X_1)\tilde N(\mathcal X)\bigr]\bigr| \leq\bbe
\Biggl[\bigl|{F}(X_1)\bigr|{\mathbf{1}}_{S_N}\sum
_{i=0}^{N-1} {e_p(X_{T_{{i}}},X_{T_{{i+1}}},
T_{{i+1}}-T_{{i}})} \Biggr],
\end{eqnarray}
{which can be bounded by ${\gamma}
E [|{F}(X_1)| ]$.}
\end{pf*}

\textit{Simulation of L\'evy bridges.} 
The adaptive method presented in this section requires fast simulation
from the bridge
distribution of {$X_{t/2}$ conditioned to $X_{t}=y$ (with $t>0$)},
whose density is given by (\ref{Eq:BrdDsity1}). We
now propose a simple yet efficient method for simulating
from the bridge distribution, valid for L\'evy processes with unimodal
density at
all times. As remarked in Section~\ref{estgen}, a sufficient condition
for a L\'evy process to have a unimodal
density for all $t>0$ is that it belongs to the class of
self-decomposable processes which includes most of the
parametric models used in the literature. The algorithm is based on
the following simple estimate.
%
\begin{prop}
Let $X$ be a L\'evy process such that the density {$f_t$} of $X_t$ is
unimodal for all $t>0$. Then,
%
%
\begin{eqnarray}
{f}^{\mathrm{br}}_{t/2}(x,y)\leq\frac{{f}_{t/2}(y/2)}{{f}_t(y)}\max\bigl\{
{f}_{t/2}(x),{f}_{t/2}(y-x)\bigr\}.\label{rejimproved}
\end{eqnarray}
\end{prop}
\begin{pf}
For all $x$ and $y$,
\begin{eqnarray*}
{f}_{t/2}(x){f}_{t/2}(y-x) = \max\bigl\{{f}_{t/2}(x),{f}_{t/2}(y-x)
\bigr\}\min\bigl\{ {f}_{t/2}(x),{f}_{t/2}(y-x)\bigr\}.
\end{eqnarray*}
By the assumption of unimodality, the density {$f_t$} may not have a
local minimum, hence, for all $a,b$, {$
\min({f}_{t/2}(a),{f}_{t/2}(b))\leq{f}_{t/2} (\frac{a+b}{2} )
$}
and the result follows.
\end{pf}

{As a consequence of the previous result,
random variates with density $f^{\mathrm{br}}_{t/2}(x,y)$ can be
simulated using the classical rejection method (Devroye \cite{devroye}), with the
proposal density given by ${\bar{f}}(x) = \frac{1}{2}({f}_{t/2}(x) +
{f}_{t/2}(y-x))$, provided that the following two requirements are met:}
\begin{itemize}[(a)]
\item[(a)] random variates with density ${f}_t(x)$ can be simulated in
bounded time;
\item[(b)] the density ${f}_t(x)$ is known explicitly or can be
evaluated in
bounded time.
\end{itemize}
Assumptions (a) and (b) are satisfied, for example, for
the variance gamma process, normal inverse Gaussian process, or for
stable processes.
Simulating a random variable
$X$ with density ${\bar{f}}(x)= \frac{1}{2}({f}_{t/2}(x) +
{f}_{t/2}(y-x))$ is straightforward: simulate a random variate
$Z$ with density ${f}_{t/2}$ and an independent Bernoulli random variate
$U$; then, take $X = Z$ if $U=0$ and $X = y-Z$ otherwise.

The expected number of iterations {needed until the}
acceptance for a {given value of} $y$ is
equal to $C = \frac{2{f}_{t/2}(y/2)}{{f}_t(y)}$. This number is bounded
for L\'evy processes with Pareto tails such as
stable. For processes with lighter tails, it may be unbounded for
large $y$, but
the probability of having a large value of $y$ in an adaptive
simulation is very small. For example, if we want to simulate
$X_{t/2}$ and $X_t$ by first simulating $X_t$ and then $X_{t/2}$ from
the bridge
law using formula \eqref{rejimproved}, we find that the conditional
expectation of the number of iterations
given $X_t$ equals $\frac{2{f}_{t/2}(X_t/2)}{{f}_t(X_t)}$, and the
unconditional expectation is
\[
\mathbb{E} \biggl[\frac{2{f}_{t/2}(X_t/2)}{{f}_t(X_t)} \biggr] = 2\int_{\mathbb R}
{f}_{t/2}(x/2)\,\mathrm{d}x = 4.
\]

\section{Numerical illustrations}\label{Sec:NumIllstr}
In this section, to simplify the discussion, we assume that the
interval $D$ is of the form $D = (-\infty,b)$.
For the numerical implementation of Algorithm \ref{ConstrSkeleton}
{given in Section~\ref{adaptive.sec},}
one needs to be able to perform the following computations efficiently:
\begin{itemize}
\item Simulation of the increments of $X_t$ for arbitrary $t$;
\item Evaluation of the density ${f_t}$ of $X_t$ for arbitrary $t$;
\item Evaluation of the ``incomplete convolution'' of the L\'evy
density: $\mathcal C(b,y):=\break \int_b^\infty s(v)s(y-v) \,\mathrm{d}v$;
\item Evaluation of the error bound $e_p(x,y,t)$, appearing in
Algorithm \ref{ConstrSkeleton}.
\end{itemize}

These computations can be performed
relatively easily, for example, for $\alpha$-stable L\'evy processes with
L\'evy density {$s(x) = |x|^{-\alpha-1}(c_- {\mathbf{1}}_{x<0} + c_+ {\mathbf{
1}}_{x>0})$
and for} the variance gamma process with L\'evy density
$s(x) = |x|^{-1}(c\mathrm{e}^{-\lambda_- |x|}{\mathbf{1}}_{x<0} + c\mathrm{e}^{-\lambda_+
|x|}{\mathbf{1}}_{x>0})$.
For $\alpha$-stable processes, {the} increments can be simulated with an
explicit algorithm {(cf. Chambers, Mallows and Stuck \cite{chambers})}, {the} density can be
computed using a
rapidly convergent series (Samorodnitsky and Taqqu \cite{taqqu})
or expressed via special
functions {(cf. G{\'o}rska and Penson \cite{gorska2011levy})}, tabulated for $t=1$ and
computed by {the} scaling property for other values of $t$. The incomplete
convolution is given by
%
%
\begin{eqnarray}
\mathcal C(b,y) = c_+ c_- b^{-1-2\alpha} B(1+2\alpha,1) F \biggl(1+\alpha,1+2
\alpha, 2+2\alpha,\frac{y}{b} \biggr),\label{stableconv}
\end{eqnarray}
where $B$ is the beta function and $F$ is the hypergeometric function,
for which a rapidly converging series is available (Gradshetyn and Ryzhik \cite{Grad}) and
which can also be tabulated prior to the Monte Carlo computation.
For the variance gamma process, the density is explicit and the
increments are straightforward to simulate (Cont and Tankov \cite{levybook}). The
incomplete convolution is given by
\[
\mathcal C(b,y) = \frac{c^2}{y} \bigl\{\mathrm{e}^{-y\lambda_+ } \operatorname{Ei}
\bigl(\lambda(b-y)\bigr) - \mathrm{e}^{y\lambda_- } \operatorname{Ei}(\lambda b)
\bigr\},
\]
where $\operatorname{Ei}(x):= \int_x^\infty\frac{\mathrm{e}^{-z}}{z}\,\mathrm{d}z$, which can
also be tabulated, and {$\lambda:=\lambda_- + \lambda_+$}.
The error bound $e_p$ for the $\alpha$-stable or the variance gamma
process can be obtained along the lines of the general computation of
{Section~\ref{estgen}} or the specific computation for the Cauchy
process in the {Appendix} \ref{cauchy.sec}.

For the numerical simulations in this section, we shall concentrate on
the Cauchy process, which is an $\alpha$-stable process with $c_+ =
c_-:= c$ and $\alpha= 1$. For this process, formula
\eqref{stableconv} simplifies to
\[
\mathcal C(b,y) = \frac{c^2}{3b^3} \Biggl\{1+3\sum
_{n=1}^\infty \frac{n+1}{n+3} \biggl(
\frac{y}{b} \biggr)^n \Biggr\} = \frac{c^2}{b^3} \biggl\{1+
\frac{b}{y} + \frac{2b^2}{y^2} + \frac{y}{b-y} + \frac{2b^3}{y^3}
\log \biggl(1-\frac{y}{b} \biggr) \biggr\}.
\]
Note that for small $y$, the series representation has more stable
behavior than the exact formula.
The error estimate $e_p$ is computed as
explained in Section~\ref{cauchy.sec} of the \hyperref[app]{Appendix}. {In both
examples below, we take $c=1$.}

\begin{example}\label{ex1}
In our first example, we evaluate the
probability
{$\bbp[\sup_{0\leq s \leq1} X_s \leq1]=\bbp(\tau>1)$, which
can be expressed in terms of the function (\ref{killed}) by taking
$T=1$, $F(X_{1})=1$, and the domain $(a,b)=(-\infty,1)$. Note that
in this case,} the starting value of the process is relatively far
from the boundary, and
hence the advantage of using the adaptive algorithm is less
important. The process will typically cross the boundary by a large
jump with a {large} overshoot, which makes the exit easy to detect,
even with a uniform discretization.

We study the performance of our adaptive algorithm for various
values of ${\gamma}$, and compare it to the standard uniform
discretization. When interpreting the results of simulations, one
needs to distinguish between the actual error (i.e., the difference
between the computed value and the true value), and the theoretical
value of the bias (computed as explained in Remark \ref{betterbias}
{above}), which
does not require the knowledge of the true value. As an estimate of
the true value, we use the value computed in an
independent simulation by uniform
discretization with 16\,384 points and $10^7$ trajectories, which is
approximately equal to $0.38935$ with a standard deviation of
$10^{-4}$. The difference between the values for 8192 and 16\,384 points
(on the
same trajectories) is smaller than $10^{-4}$, hence one can
presume that, {for all practical purposes,} convergence up to this
precision has been achieved.

Figure~\ref{actual.fig} shows the dependence of the values computed by
the two algorithms on the computational time required for $10^6$
MC trajectories, for different numbers of discretization points (for
the uniform discretization) and different values of the tolerance
parameter ${\gamma}$ (for the adaptive algorithm). While the
uniform discretization algorithm exibits a clear bias which decreases
as the number of discretization dates increases, the adaptive
algorithm removes the bias completely; all values returned by this
algorithm are within confidence bounds of the true value.

\begin{figure}

\includegraphics{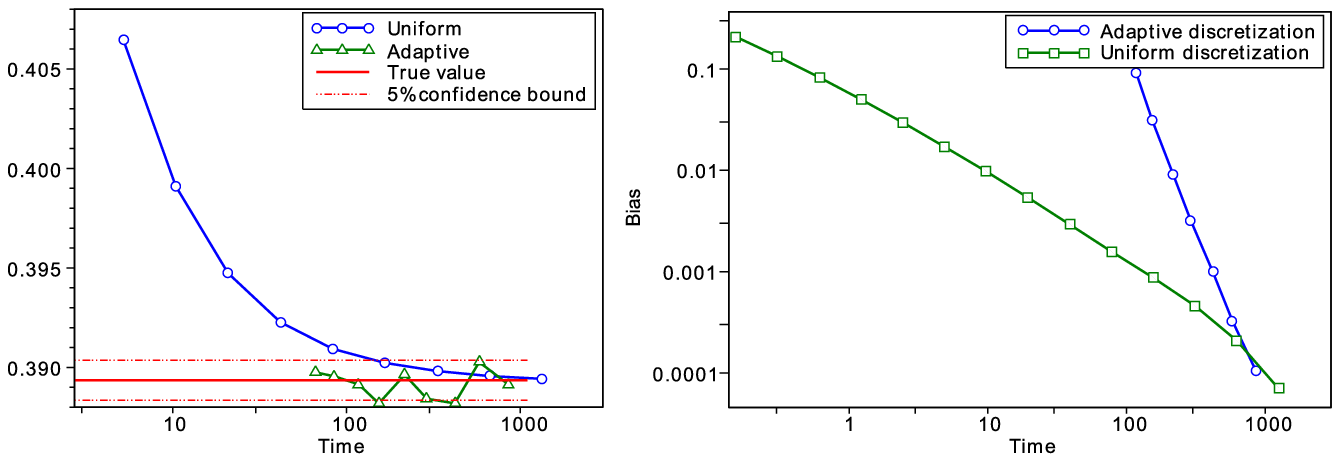}

\caption{Illustration for Example \protect\ref{ex1}. Left: values returned by the
uniform discretization algorithm
and the adaptive algorithm, as function of the computational time for
$10^6$ paths,
measured in seconds. Different points on the graph correspond
to different numbers of discretization dates for the uniform
discretization (ranging from 32 to 8192) and different values of the
tolerance parameter {$\gamma$} for the adaptive algorithm (ranging
from $7$ to $7\times1 0^{-4}$). The curve for the uniform
discretization is smooth because all the points have been generated
using the same trajectories, while for the adaptive discretization
different paths have been used. Right: comparison of the theoretical
bias of the adaptive algorithm with the actual discretization bias of
the uniform discretization.}
\label{actual.fig}
\end{figure}

The theoretical bias, computed as explained in Remark
\ref{betterbias}, is greater than the actual error, because the error
estimates of {\hyperref[app]{Appendix}} \ref{cauchy.sec} are upper bounds, and because it
does not take into account the possible cancellation of errors on
different intervals.
Figure~\ref{actual.fig}, right graph, compares the theoretical
estimate of the bias of
the adaptive algorithm with the actual bias of the uniform
discretization. One can see that for small computational times, the
theoretical bias for the adaptive algorithm is greater than the error
of the uniform discretization, however, the theoretical bias converges
to zero much faster, and for relatively large computational times is
actually smaller than the error of the uniform discretization. The
empirical convergence rate (estimated from the slope
of the straight lines) is $T^{-0.81}$ for the uniform discretization
and $T^{-3.4}$ for the theoretical bias of the adaptive algorithm.
\end{example}

\begin{figure}

\includegraphics{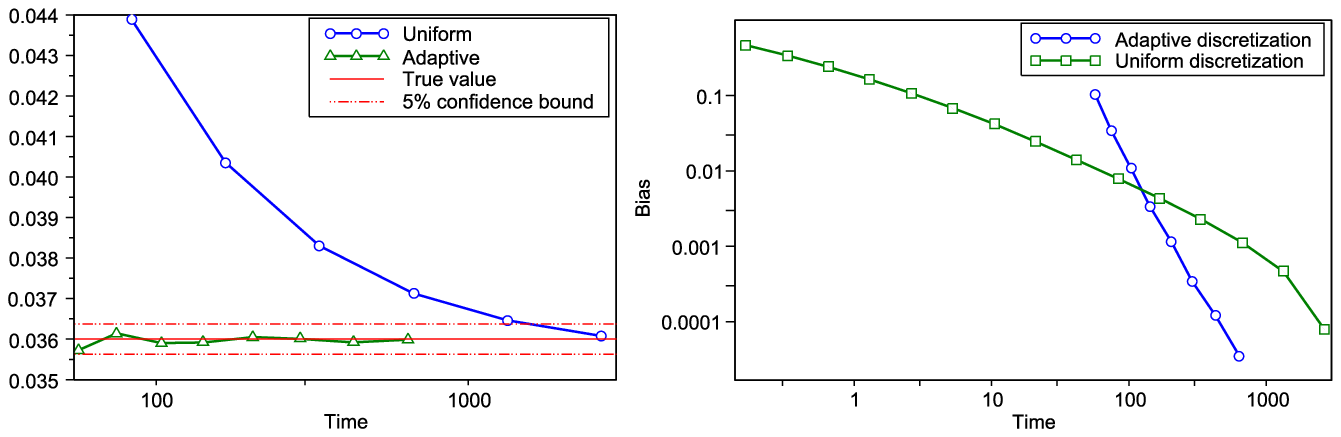}

\caption{Illustration for Example \protect\ref{ex2}. Left: values returned by the
uniform discretization algorithm
and the adaptive algorithm, as function of the computational time for
$10^6$ paths,
measured in seconds. Different points on the graph correspond
to different numbers of discretization dates for the uniform
discretization (ranging from 256 to 16384) and different values of the
tolerance parameter {$\gamma$} for the adaptive algorithm (ranging
from $9$ to $9\times10^{-3}$). Right: comparison of the theoretical
bias of the adaptive algorithm with the actual discretization bias of
the uniform discretization.}
\label{actual2.fig}
\end{figure}

\begin{example}\label{ex2}
In our second example, we evaluate the probability
{$\bbp[\sup_{0\leq s \leq1} X_s \leq10^{-2}]$, which again can be
expressed in terms of the function (\ref{killed}) by taking $T=1$,
$F(X_{1})=1$, and the domain $(a,b)=(-\infty,10^{-2})$. In contrast to
Example \ref{ex1}, here} we consider a
situation where the starting point is close to the
boundary. In this case, as
we shall see below, the advantage of the adaptive algorithm is
more {striking}, since the process can cross the boundary and come back
while it is still close to the starting point and, hence, a very fine
discretization will be necessary to detect this event with uniformly
spaced observations. As a result, for the uniform discretization we do
not observe convergence to
a sufficient precision even with 16\,384 points, and therefore the
true value cannot be estimated as in the previous example. Instead, we
shall use as the true value the value produced by the adaptive
algorithm with $10^7$ Monte Carlo paths and equal to $0.0360$, with
standard deviation of $6\times10^{-5}$ and theoretical bias of
$3\times10^{-5}$.

Similarly to the previous example, Figure~\ref{actual2.fig} shows the
dependence of the values computed by
the two algorithms on the computational time required for $10^6$
MC trajectories. Here, the adaptive algorithm exhibits the same kind
of behavior {as in the Example \ref{ex1} above}: all the points generated by
the algorithm are within the
confidence bounds of the true value. However, for the uniform
discretization, the convergence is much slower than before and only
the last value obtained with 16384 discretization points falls within
the confidence bounds. Figure~\ref{actual2.fig}, right graph, compares
the theoretical
estimate of the bias of
the adaptive algorithm with the actual bias of the uniform
discretization. Once again, the behavior of the adaptive algorithm is
roughly the same as in the previous example, showing that the method
is robust with respect to the parameters on the problem. On the other
hand, as expected, the uniform
discretization presents a significant bias in this case (the
convergence rates are similar to those obtained in the previous
example, but the constant for the uniform discretization is much bigger).
\end{example}

\begin{appendix}\label{app}

\section{Proofs of Section \texorpdfstring{\protect\ref{asymp}}{2}}\label{Sec:ProofTh:STBS}

\subsection{Proof of Theorem \texorpdfstring{\protect\ref{Th:STBS}}{2.1}}\label{SbSec:PrMR}

{Throughout the proof, we shall use the {notation}
%
%
\begin{equation}
\label{StndNotMaxMin} {\bar{Y}_{t}:=\sup_{0\leq s\leq t}
Y_{s} \quad\mbox{and}\quad \underline{Y}_{t}:=\inf
_{0\leq s\leq t} Y_{s}}
\end{equation}
for a given c\'adl\'ag process $(Y_{t})_{t\geq 0}$. Without} loss of
generality (by considering
separately the positive and the negative part), we can and will assume
that $\varphi$ is nonnegative. Additionally, assume that $a\in(-\infty
,0)$ and $b\in(0,\infty)$. The
cases $a=-\infty$ and $b=\infty$ will be evident from the
proof below. We also let
{$\|\varphi\|_{\infty}:=\operatorname{ess\, sup}_{x}\varphi(x)$}, {$\|\varphi
\|_\mathrm{ Lip}$ be the Lipschitz norm\vspace*{1pt} of $\varphi$},
$I_{\delta}(y):=({y-\delta},y+\delta)$, $\eta:=\delta_{0}/2$,
$c=b\wedge|a|$, $B:=\{\tau\leq t\}=\{\bar{X}_{t}\geq b \mbox{ or
}\underline{X}_{t}\leq a\}$, $U_{t}^{\varepsilon}:=\sup_{s\leq
t}|X^\varepsilon_{s}|$, and {$a_{\varepsilon}:=\sup_{x} s_{\varepsilon
}(x)$, which} are finite in light of (\ref{Cond1}). {In what follows,
$\calF_{t}^{\varepsilon}:=\sigma (X^{\varepsilon}_{s}\dvt s\leq t
)\vee\calN$ where $\calN$ denotes the null sets of $\calF$. {To
lighten the notation below, whenever the $\operatorname{ess\, sup}$ of a
function $g$, defined $\calL$-a.e. in some region, is considered, we
shall simply write $\sup_{u}g(u)$ instead of $\operatorname{ess\,
sup}_{u}g(u)$.}}\vadjust{\goodbreak}

{The idea is to condition on the number of jumps of the
compound Poisson component $Z^{\varepsilon}$. To this end, let us denote
\begin{eqnarray*}
A_{k}(t)&=& \bbe \bigl( \varphi(X_{\tau}){
\mathbf{1}}_{\{\tau\leq t, X_{t}\in
I_{\delta}(y),N_{t}^{\varepsilon}=k\}} \bigr)\qquad\mbox{for }k=0,1,2,
\\
A_{3}(t)&=& \bbe \bigl( \varphi(X_{\tau}){
\mathbf{1}}_{\{\tau\leq t, X_{t}\in
I_{\delta}(y),N_{t}^{\varepsilon}\geq3\}} \bigr),
\end{eqnarray*}
so that clearly
%
%
\begin{equation}
\label{T1-T4} \bbe\bigl(\varphi(X_{\tau})\mathbf{1}_{\{\tau\leq t, X_{t}\in I_{\delta}(y)\}
}
\bigr)=A_{0}(t)+\cdots+A_{3}(t).
\end{equation}
Note that each of the terms on the right-hand side of the previous
{equation} can be expressed as
%
%
\begin{equation}
\label{WAkIPk} A_{k}(t)=\int_{y-\delta}^{y+\delta}
P_{t}^{k}(u) \,\mathrm{d}u \qquad(k=0,\ldots,3)
\end{equation}
for some nonnegative functions $P^{k}_{t}(u)$. Indeed, for $k=0,1,2$,
by the standard definition {of conditional} expectation,
%
%
\begin{eqnarray}\label{WAkIPkb}
\nonumber
A_{k}(t)&=& \bbe \bigl( \varphi(X_{\tau}){
\mathbf{1}}_{\{\tau\leq t, X_{t}\in
I_{\delta}(y),N_{t}^{\varepsilon}=k\}} \bigr)=\bbe \bigl( \bbe \bigl(\varphi(X_{\tau}){
\mathbf{1}}_{\{\tau\leq t, N_{t}^{\varepsilon}=k\}}\vert X_{t} \bigr)\mathbf{1}_{\{X_{t}\in I_{\delta}(y)\}} \bigr)
\nonumber
\\[-8pt]
\\[-8pt]
\nonumber
&=&\int_{{y-\delta}}^{y+\delta}\bbe\bigl(
\varphi(X_{\tau})\mathbf{1}_{\{
\tau\leq t,N_{t}^{\varepsilon}=k\}}\vert X_{t}=u
\bigr)f_{t}(u)\,\mathrm{d}u=: \int_{{y-\delta}}^{y+\delta}
P_{t}^{k}(u) \,\mathrm{d}u.
\end{eqnarray}
The case $k=3$ is treated in the same way.
Let} us analyze each of the four terms in the right-hand side of (\ref{T1-T4}).

(1) \textit{No big jump}.
Note that, {on the event} $N_{t}^{\varepsilon}=0$,
$X_{s}=X_{s}^{\varepsilon}$ for all $s\leq t$ and, {thus, $\{\tau\leq
t\}=\{\tau^{\varepsilon}\leq t\}$, where $\tau^{\varepsilon}:=\inf\{
u\geq 0\dvt X_{u}^{\varepsilon}\notin(a,b)\}$. Therefore,
\begin{eqnarray*}
A_{0}(t)&=&\bbe \bigl( \varphi\bigl(X_{\tau^{\varepsilon}}^{\varepsilon}
\bigr)\mathbf{ 1}_{\{\tau^{\varepsilon}\leq t, X^{\varepsilon}_{t}\in I_{\delta
}(y),N_{t}^{\varepsilon}=0\}} \bigr)\\
&=&\bbe \bigl( \varphi\bigl(X_{\tau
^{\varepsilon}}^{\varepsilon}
\bigr)\mathbf{1}_{\{\tau^{\varepsilon}\leq t,
X^{\varepsilon}_{t}\in I_{\delta}(y)\}} \bigr)\bbp \bigl(N_{t}^{\varepsilon}=0
\bigr),
\end{eqnarray*}
where in the last equality we used the independence of $X^{\varepsilon
}$ and {$N^{\varepsilon}$}.}
Next, conditioning on $\calF_{\tau^{\varepsilon}}^{\varepsilon}$, it
follows that
\begin{eqnarray*}
A_{0}(t)= \mathrm{e}^{-\lambda_{\varepsilon}t}\bbe\bigl(\varphi\bigl(X_{\tau^{\varepsilon
}}^{\varepsilon}
\bigr)\mathbf{1}_{\{\tau^{\varepsilon}\leq t,
X_{t}^{\varepsilon}\in I_{\delta}(y)\}}\bigr)= \mathrm{e}^{-\lambda_{\varepsilon}t} \bbe \bigl( \bbe\bigl(
\mathbf{1}_{\{X_{t}^{\varepsilon}\in I_{\delta
}(y)\}}\vert \calF_{\tau^{\varepsilon}}^{\varepsilon} \bigr)
\varphi \bigl(X_{\tau^{\varepsilon}}^{\varepsilon}\bigr)\mathbf{1}_{\{\tau^{\varepsilon}\leq
t\}} \bigr).
\end{eqnarray*}
By Markov's property,
\begin{eqnarray*}
A_{0}(t)&=& \mathrm{e}^{-\lambda_{\varepsilon}t} \bbe \bigl(\bbe \bigl({
\mathbf{1}}_{\{X_{t}^{\varepsilon}-X_{\tau
^{\varepsilon}}^{\varepsilon}+X_{\tau^{\varepsilon}}^{\varepsilon}\in
I_{\delta}(y)\}}\vert \calF^{\varepsilon}_{\tau^{\varepsilon}} \bigr)\varphi
\bigl(X_{\tau^{\varepsilon}}^{\varepsilon}\bigr)\mathbf{1}_{\{\tau
^{\varepsilon}\leq t\}} \bigr)
\\
&=& { \mathrm{e}^{-\lambda_{\varepsilon}t}} \bbe \bigl(F \bigl(X_{\tau^{\varepsilon}}^{\varepsilon},t-
\tau ^{\varepsilon} \bigr)\varphi\bigl(X_{\tau^{\varepsilon}}^{\varepsilon}\bigr)\mathbf{
1}_{\{\tau^{\varepsilon}\leq t\}} \bigr),
\end{eqnarray*}
where $F(z,s)=\bbp (z+X_{s}^{\varepsilon}\in I_{\delta}(y) )$.
{Note that if $\tau^{\varepsilon}=t$, then $F (X_{\tau^{\varepsilon
}}^{\varepsilon},t-\tau^{\varepsilon} )=0$ since $X_{\tau
^{\varepsilon}}^{\varepsilon}\in(a,b)^{c}$ and $I_{\delta}(y)\subset
(a,b)$. On the other hand, on the event $\tau^{\varepsilon}<t$,
\begin{eqnarray*}
F \bigl(X_{\tau^{\varepsilon}}^{\varepsilon},t-\tau^{\varepsilon} \bigr)=\int
_{{y-\delta}}^{y+\delta} f_{t-\tau^{\varepsilon}}^{\varepsilon
}
\bigl(u-X_{\tau^{\varepsilon}}^{\varepsilon}\bigr) \,\mathrm{d}u\leq\int_{y-\delta
}^{y+\delta}
\sup_{0< s \leq t} \sup_{x\in(a,b)^c} f_{s}^{\varepsilon}(u-x)\,\mathrm{d}u,
\end{eqnarray*}
since again $X_{\tau^{\varepsilon}}^{\varepsilon}\in(a,b)^{c}$.
Putting the two previous cases together and recalling (\ref{WAkIPk}),
we have
%
%
\begin{eqnarray}\label{Eq:FrsTrmRmd}
\nonumber
A_{0}(t)= \int_{y-\delta}^{y+\delta}
P_{t}^{0}(u)\,\mathrm{d}u&\leq &\int_{y-\delta}^{y+\delta}
\Bigl(\mathrm{e}^{-\lambda_{\varepsilon}t}\|\varphi\| _{\infty} \sup_{0< s \leq t}
\sup_{x\in(a,b)^c} f_{s}^{\varepsilon
}(u-x) \Bigr)\,\mathrm{d}u
\nonumber
\\[-8pt]
\\[-8pt]
\nonumber
&=:&\int_{y-\delta}^{y+\delta} \bar{P}^{0}_{t}(u)\,\mathrm{d}u,
\end{eqnarray}
implying that $P^{0}_{t}(u)\leq\bar{P}_{t}^{0}(u)$, for $\calL$-a.e.
$u\in(a+\delta_{0},b-\delta_{0})$. Furthermore,
using (\ref{KyInBn})(ii),}
\[
\sup_{{a+\delta_{0}<u< b-\delta_{0}}}P_{t}^{0}(u)\leq \sup
_{{a+\delta_{0}<u< b-\delta_{0}}}\bar{P}_{t}^{0}(u)\leq\|\varphi\|
_{\infty}c_{3}(\delta_0,\varepsilon)
t^{3}\qquad (t<t_{0}).
\]

(2) \textit{One big jump}.
Let $\tau_{i}$ and $Y_{i}$ be the time and size of the $i${th} jump of
$Z^{\varepsilon}$.
Clearly, on the event $\{N_{t}^{\varepsilon}=1\}$, 
\begin{eqnarray*}
\varphi(X_{\tau})\mathbf{1}_{\{\tau\leq t, X_{t}\in I_{\delta
}(y),{N_{t}^{\varepsilon}=1}\}} &=&\varphi\bigl(X_{\tau}^{\varepsilon}
\bigr)\mathbf{1}_{\{\tau< \tau_{1},
X_{t}^{\varepsilon}+Y_{1}\in I_{\delta}(y),{N_{t}^{\varepsilon}=1}\}}
\\
&&{} +\varphi\bigl(X_{\tau}^{\varepsilon}+Y_{1}\bigr){
\mathbf{1}}_{\{\tau_{1}\leq
\tau\leq t, X_{t}^{\varepsilon}+Y_{1}\in I_{\delta
}(y),{N_{t}^{\varepsilon}=1}\}}
\\
&\leq&\|\varphi\|_{\infty}\mathbf{1}_{\{X_{t}^{\varepsilon}+Y_{1}\in
I_{\delta}(y),{N_{t}^{\varepsilon}=1}\}} \mathbf{1}_{\{\bar{X}^{\varepsilon}_{t}\geq b\ \mathrm{or}\ \underline
{X}_{t}^{\varepsilon}\leq a\}}
\\
&&{} +\|\varphi\|_{\infty}\mathbf{1}_{\{X_{t}^{\varepsilon}+Y_{1}\in
I_{\delta}(y),{N_{t}^{\varepsilon}=1}\}} {\mathbf{1}_{\{\bar{X}^{\varepsilon}_{t}+Y_{1}\geq b\ \mathrm{or}\
\underline{X}_{t}^{\varepsilon}+Y_{1}\leq a\}}}.
\end{eqnarray*}
It follows that
\begin{eqnarray*}
0\leq A_{1}(t) &\leq&\|\varphi\|_{\infty} \bbe ({
\mathbf{1}}_{ \{U_{t}^{\varepsilon}\geq c,
X_{t}^{\varepsilon}+Y_{1}\in I_{\delta}(y),N_{t}^{\varepsilon}=1 \}
} )
\\
&&{} + \|\varphi\|_{\infty} \bbe (\mathbf{1}_{ \{Y_{1}\geq b-\bar{X}_{t}^{\varepsilon}\ \mathrm{or}\ Y_{1}\leq a-\underline{X}^{\varepsilon}_{t} \}}{
\mathbf{1}}_{ \{
X_{t}^{\varepsilon}+Y_{1}\in I_{\delta}(y),N_{t}^{\varepsilon}=1 \}
} )
\\
&=& \underbrace{\mathrm{e}^{-\lambda_{\varepsilon}t}\|\varphi\|_{\infty} \lambda
_{\varepsilon}t \bbe (\mathbf{1}_{ \{U_{t}^{\varepsilon}\geq c,
X_{t}^{\varepsilon}+Y_{1}\in I_{\delta}(y) \}} )}_{A_{1,1}(t)}
\\
&&{} + \underbrace{{\mathrm{e}^{-\lambda_{\varepsilon}t}\|\varphi\|_{\infty}\lambda
_{\varepsilon}t \bbe (\mathbf{1}_{ \{Y_{1}\geq b-\bar{X}_{t}^{\varepsilon}\ \mathrm{or}\ Y_{1}\leq a-\underline{X}^{\varepsilon}_{t} \}}\mathbf{1}_{ \{
X_{t}^{\varepsilon}+Y_{1}\in I_{\delta}(y) \}}
)}}_{A_{1,2}(t)},
\end{eqnarray*}
{where in the last equality we use the joint independence of
$N^{\varepsilon}$, $Y_{1}$, and $X^{\varepsilon}$.} Conditioning on
$\sigma(X_{s}^{\varepsilon}\dvt s\geq 0)$ {and applying Fubini},
%
%
\begin{eqnarray}\label{Eq:DefnP11}
\nonumber
A_{1,1}(t)&=&{\mathrm{e}^{-\lambda_{\varepsilon}t}}\|\varphi\|_{\infty}t \bbe
\biggl(\mathbf{1}_{\{{U}_{t}^{\varepsilon}\geq c\}} \int_{{y-\delta}-X_{t}^{\varepsilon}}^{y+\delta-X_{t}^{\varepsilon
}}s_{\varepsilon}(v)\,\mathrm{d}v
\biggr)
\nonumber
\\[-8pt]
\\[-8pt]
\nonumber
&=&\int_{{y-\delta}}^{y+\delta}\underbrace{{\mathrm{e}^{-\lambda_{\varepsilon
}t}}\|
\varphi\|_{\infty}t \bbe \bigl(\mathbf{1}_{\{{U}_{t}^{\varepsilon}\geq
c\}}s_{\varepsilon}
\bigl(u-X_{t}^{\varepsilon}\bigr) \bigr)}_{\bar
{P}_{t}^{1,1}(u)}\,\mathrm{d}u.
\end{eqnarray}
Using (\ref{Cond1}) and Lemma \ref{KyLm1},
\(
\sup_{u}\bar{P}_{t}^{1,1}(u)\leq{\mathrm{e}^{-\lambda_{\varepsilon}t}}t \|
\varphi\|_{\infty}a_{\varepsilon}\bbp ({U}_{t}^{\varepsilon}\geq
c )\leq \mathrm{e}^{-\lambda_{\varepsilon}t}
a_{\varepsilon} \|\varphi\|_{\infty}\times  C_{2}(c,\varepsilon) t^{3},
\)
where $\varepsilon>0$ is chosen small enough. Similarly, {conditioning
on $\sigma(X_{s}^{\varepsilon}\dvt s\geq 0)$, making the substitution
$u=X_{t}^{\varepsilon}+v$, and applying Fubini,}
%
%
\begin{eqnarray}\label{Eq:DefnP12}
A_{1,2}(t) &=&{\mathrm{e}^{-\lambda_{\varepsilon}t}}\|\varphi\|_{\infty}t\bbe \biggl(
\int\mathbf{ 1}_{\{v\leq a-\underline{X}^{\varepsilon}_{t}\ \mathrm{or}\ v\geq b-\bar
{X}_{t}^{\varepsilon}\}}\mathbf{1}_{\{{y-\delta}<  X_{t}^{\varepsilon
}+v\leq y+\delta\}}s_{\varepsilon}(v)\,\mathrm{d}v \biggr)
\nonumber
\\[-8pt]
\\[-8pt]
\nonumber
&=&\int_{{y-\delta}}^{y+\delta} \underbrace{{\mathrm{e}^{-\lambda_{\varepsilon}t}}\|
\varphi\|_{\infty}t\bbe \bigl( \mathbf{ 1}_{\{u\leq a+X_{t}^{\varepsilon}-\underline{X}^{\varepsilon}_{t}\ \mathrm{or}\ u\geq b+X_{t}^{\varepsilon}-\bar{X}_{t}^{\varepsilon}\}
}s_{\varepsilon}
\bigl(u-X_{t}^{\varepsilon}\bigr) \bigr)}_{\bar
{P}_{t}^{1,2}(u)}\,\mathrm{d}u.
\end{eqnarray}
Using again Lemma \ref{KyLm1},
%
%
\begin{eqnarray}\label{In:RepUppBnd}
\nonumber
\sup_{u\in(a+\delta_{0},b-\delta_{0})} \bar{P}_{t}^{1,2}(u)&
\leq&{\mathrm{e}^{-\lambda_{\varepsilon}t}}\|\varphi\| _{\infty} t a_{\varepsilon} \bbp
\bigl({X}_{t}^{\varepsilon}-\underline {X}_{t}^{\varepsilon}
\geq \delta_{0}\ \mathrm{or}\ \bar{X}_{t}^{\varepsilon
}-X_{t}^{\varepsilon}
\geq \delta_{0} \bigr)
\\
&\leq&{\mathrm{e}^{-\lambda_{\varepsilon}t}}\|\varphi\|_{\infty} t a_{\varepsilon}
\bbp \bigl(\bar{X}_{t}^{\varepsilon}-\underline {X}_{t}^{\varepsilon}
\geq \delta_{0} \bigr)
\nonumber
\\[-8pt]
\\[-8pt]
\nonumber
&\leq&{\mathrm{e}^{-\lambda_{\varepsilon}t}}\|\varphi\|_{\infty} t a_{\varepsilon}\bbp \Bigl(
\sup_{s\leq t}\bigl|{X}_{s}^{\varepsilon}\bigr|\geq
\delta_{0}/2 \Bigr)
\\
& \leq&{\mathrm{e}^{-\lambda_{\varepsilon}t}} \|\varphi\|_{\infty} a_{\varepsilon
}C_{2}(
\delta_0/2,\varepsilon) t^{3}.\nonumber
\end{eqnarray}
Therefore, recalling from (\ref{WAkIPk}), the nonnegative
function $P_{t}^{1}(u)$ is such that,for $\calL$-a.e. $u\in(a+\delta
_{0}, b-\delta_{0})$,
$0\leq P_{t}^{1}(u)\leq\sum_{\ell=1}^{2}\bar{P}_{t}^{1,\ell}(u)\leq\|
\varphi\|_{\infty}a_{\varepsilon}t^{3} ( C_{2}(c,\varepsilon
)+C_{2}(\eta,\varepsilon)  )$.

(3) \textit{Two big jumps.} As before, let $\tau_{i}$ and
$Y_{i}$ be the time and size of the $i${th} jump of $Z^{\varepsilon}$.
Clearly,
\begin{eqnarray*}
\varphi(X_{\tau})\mathbf{1}_{\{\tau\leq t, X_{t}\in I_{\delta
}(y),{N_{t}^{\varepsilon}=2}\}} &=&\varphi\bigl(X_{\tau}^{\varepsilon}
\bigr)\mathbf{1}_{\{\tau< \tau_{1},
X_{t}^{\varepsilon}+Y_{1}+Y_{2}\in I_{\delta}(y),{N_{t}^{\varepsilon
}=2}\}}
\\
&&{} +\varphi\bigl(X_{\tau}^{\varepsilon}+Y_{1}\bigr){
\mathbf{1}}_{\{\tau_{1}\leq
\tau<\tau_{2}, X_{t}^{\varepsilon}+Y_{1}+Y_{2}\in I_{\delta
}(y),{N_{t}^{\varepsilon}=2}\}}
\\
&&{} +\varphi\bigl(X_{\tau}^{\varepsilon}+Y_{1}+Y_{2}
\bigr)\mathbf{ 1}_{\{\tau_{2}\leq \tau\leq t,
X_{t}^{\varepsilon}+Y_{1}+Y_{2}\in I_{\delta}(y),{N_{t}^{\varepsilon
}=2}\}}
\\
& \leq& {\|\varphi\|_{\infty}} {\mathbf{1}}_{\{\exists s<\tau_1\dvt X^\varepsilon_s
\notin(a,b); X_{t}^{\varepsilon}+Y_{1}+Y_{2}\in I_{\delta
}(y);{N_{t}^{\varepsilon}=2}\}}
\\
&&{} + {\varphi\bigl(X^\varepsilon_\tau+ Y_1\bigr) {{
\mathbf{1}}}_{\{\exists s\in
[\tau_1,\tau_2)\dvt X^\varepsilon_s + Y_1
\notin(a,b); X_{t}^{\varepsilon}+Y_{1}+Y_{2}\in I_{\delta
}(y);{N_{t}^{\varepsilon}=2}\}}}
\\
&&{} + {\|\varphi\|_{\infty}} {\mathbf{1}}_{\{\exists s\in
[\tau_2,t]\dvt X^\varepsilon_s + Y_1+Y_2
\notin(a,b); X_{t}^{\varepsilon}+Y_{1}+Y_{2}\in I_{\delta
}(y);{N_{t}^{\varepsilon}=2}\}}.
\end{eqnarray*}
{Then, using the independence of $N^{\varepsilon}$, the $Y_{i}$'s, and
$X^{\varepsilon}$ in the first and last terms, we have the inequality}:
%
%
\begin{eqnarray}\label{Eq:UpBnd1}
A_{2}(t) &\leq& {\mathrm{e}^{-\lambda_{\varepsilon}t}}\bigl(t^{2}/2\bigr)
\lambda_{\varepsilon}^{2}\| \varphi\|_{\infty}\bbe ({
\mathbf{1}}_{\{{U}_{t}^{\varepsilon}\geq
c,X_{t}^{\varepsilon}+Y_{1}+Y_{2}\in I_{\delta}(y)\}} )
\nonumber
\\
& &{}+ \bbe \bigl(\varphi\bigl( {X}_{\tau}^{\varepsilon}+Y_{1}
\bigr)\mathbf{1}_{\{\bar{X}^{\varepsilon
}_{t}+Y_{1}\geq b\ \mathrm{or}\ \underline{X}_{t}^{\varepsilon}+Y_{1}\leq
a;X_{t}^{\varepsilon}+Y_{1}+Y_{2}\in I_{\delta}(y);{N_{t}^{\varepsilon
}=2}\}} \bigr)
\nonumber
\\[-8pt]
\\[-8pt]
\nonumber
&&{} + {\mathrm{e}^{-\lambda_{\varepsilon}t}} {\bigl(t^{2}/2\bigr)\lambda_{\varepsilon}^{2}
\| \varphi\|_{\infty}\bbe (\mathbf{ 1}_{\{\bar{X}^{\varepsilon}_{t}+Y_{1}+Y_{2}\geq b\ \mathrm{or}\ \underline
{X}^{\varepsilon}_{t}+Y_{1}+Y_{2}\leq a; X_{t}^{\varepsilon
}+Y_{1}+Y_{2}\in I_{\delta}(y)\}} )},\qquad
\\
&=:& {A_{2,1}(t)+A_{2,2}(t)+A_{2,3}(t).}
\nonumber
\end{eqnarray}
As before, conditioning on $\sigma(X_{s}^{\varepsilon}\dvt s\geq 0)$,
{changing variable from $w$ to $u=X_{t}^{\varepsilon}+v+w$, and
applying Fubini},
%
%
\begin{eqnarray}\label{Eq:DfnP21}
A_{2,1}(t) &=& {\mathrm{e}^{-\lambda_{\varepsilon}t}}2^{-1}\|\varphi
\|_{\infty}t^{2}\bbe \biggl(\iint\mathbf{1}_{\{U_{t}^{\varepsilon}\geq c\}} {
\mathbf{1}}_{\{{y-\delta}<X_{t}^{\varepsilon}+w+v<y+\delta\}}s_{\varepsilon
}(v)s_{\varepsilon}(w)\,\mathrm{d}v \,\mathrm{d}w \biggr)
\nonumber
\\
&=&\int_{{y-\delta}}^{y+\delta} {\mathrm{e}^{-\lambda_{\varepsilon}t}}2^{-1}
\|\varphi\|_{\infty}t^{2} \int_{-\infty}^{\infty}s_{\varepsilon}(v)
\bbe \bigl(\mathbf{1}_{\{
U^{\varepsilon}_{t}\geq c\}}s_{\varepsilon}\bigl(u-X_{t}^{\varepsilon
}-v
\bigr) \bigr)\,\mathrm{d}v \,\mathrm{d}u
\\
&=:&{\int_{{y-\delta}}^{y+\delta} \bar{P}_{t}^{2,1}(u)\,\mathrm{d}u},\nonumber
\end{eqnarray}
and, hence,
\begin{eqnarray*}
\sup_{u}\bar{P}_{t}^{2,1}(u)&
\leq{\mathrm{e}^{-\lambda_{\varepsilon}t}}2^{-1}\| \varphi\|_{\infty}t^{2}
\lambda_{\varepsilon}a_{\varepsilon} \bbp \bigl(U_{t}^{\varepsilon}
\geq c \bigr)\leq{\mathrm{e}^{-\lambda_{\varepsilon}t}} 2^{-1}\|\varphi\|_{\infty}
\lambda_{\varepsilon}a_{\varepsilon} C_{1}(c,\varepsilon)
t^{3}.
\end{eqnarray*}
Similarly, $A_{2,3}(t)$ can be written as
%
%
\begin{eqnarray}\label{Eq:DfnP23}
&&{\mathrm{e}^{-\lambda_{\varepsilon}t}} 2^{-1}\|\varphi\|_{\infty}t^{2}
\bbe \biggl(\iint\mathbf{ 1}_{\{\bar{X}^{\varepsilon}_{t}+v+w\geq b\ \mathrm{or}\ \underline
{X}^{\varepsilon}_{t}+v+w\leq a\}} \mathbf{1}_{\{{y-\delta}<X_{t}^{\varepsilon}+w+v<y+\delta\}}s_{\varepsilon
}(v)s_{\varepsilon}(w)\,\mathrm{d}v
\,\mathrm{d}w \biggr)
\nonumber
\\
&&\quad= \int_{{y-\delta}}^{y+\delta}{\mathrm{e}^{-\lambda_{\varepsilon}t}}2^{-1}
\| \varphi\|_{\infty}t^{2}
\nonumber
\\[-8pt]
\\[-8pt]
\nonumber
&&\hspace*{22pt}\qquad{}\times\int\bbe \bigl(\mathbf{
1}_{\{\bar{X}^{\varepsilon}_{t}-X_{t}^{\varepsilon}+u\geq b\ \mathrm{or}\ \underline{X}^{\varepsilon}_{t}-X_{t}^{\varepsilon}+u\leq a\}
}s_{\varepsilon}\bigl(u-X_{t}^{\varepsilon}-v
\bigr) \bigr) s_{\varepsilon}(v)\,\mathrm{d}v \,\mathrm{d}u
\\
&&\quad=: {\int_{{y-\delta}}^{y+\delta}\bar{P}_{t}^{2,3}(u)\,\mathrm{d}u,\nonumber}
\end{eqnarray}
and, thus, as in (\ref{In:RepUppBnd}),
\begin{eqnarray*}
\sup_{u\in[a+\delta_{0},b-\delta_{0}]}\bar{P}_{t}^{2,3}(u)&\leq&
{\mathrm{e}^{-\lambda_{\varepsilon}t}}2^{-1}\|\varphi\|_{\infty}t^{2}
\lambda _{\varepsilon}a_{\varepsilon}\bbp \bigl({X}_{t}^{\varepsilon}-
\underline {X}_{t}^{\varepsilon}\geq \delta_{0}\ \mathrm{or}\
\bar{X}_{t}^{\varepsilon
}-X_{t}^{\varepsilon}\geq
\delta_{0} \bigr)
\\
&\leq&{\mathrm{e}^{-\lambda_{\varepsilon}t}} 2^{-1} \|\varphi\|_{\infty}
\lambda_{\varepsilon}a_{\varepsilon} C_{1}(\delta_0/2,
\varepsilon) t^{3}.
\end{eqnarray*}
Finally, we provide an upper bound for $A_{2,2}(t)$. First, we use the
bound $\varphi(X_{\tau}^{\varepsilon}+Y_{1})\leq\varphi(Y_{1})+\|
\varphi\|_{\mathrm{ Lip}}U_{t}^{\varepsilon}$ and again the independence of
$N^{\varepsilon}$, the $Y_{i}$'s, and $X^{\varepsilon}$ to get
\[
A_{2,2}(t)\leq \mathrm{e}^{-\lambda_{\varepsilon}t}\bigl(t^{2}/2\bigr)
\lambda_{\varepsilon
}^{2}\bbe \bigl( \bigl\{\varphi(Y_{1})+
\|\varphi\|_{\mathrm{
Lip}}U_{t}^{\varepsilon} \bigr\}{
\mathbf{1}}_{\{\bar{X}^{\varepsilon
}_{t}+Y_{1}\geq b\ \mathrm{or}\ \underline{X}_{t}^{\varepsilon}+Y_{1}\leq
a;X_{t}^{\varepsilon}+Y_{1}+Y_{2}\in I_{\delta}(y)\}} \bigr).
\]
Next, by conditioning on
$\sigma(X_{s}^{\varepsilon}\dvt s\geq 0)\vee\sigma(Y_{1})$, we may
write\footnote{Here and below we use
the convention $(x,y) = \varnothing$ {and $(x,y)^{c}=(-\infty,\infty)$
for $x>y$}.} $A_{2,2}(t)$ as
\begin{eqnarray*}
&& {\mathrm{e}^{-\lambda_{\varepsilon}t}}\bigl(t^{2}/2\bigr) \lambda_{\varepsilon}\bbe
\biggl( \bigl\{\varphi(Y_{1})+\|\varphi\|_{\mathrm{ Lip}}U_{t}^{\varepsilon}
\bigr\}\mathbf{1}_{\{\bar{X}^{\varepsilon}_{t}+Y_{1}\geq b\ \mathrm{or}\ \underline
{X}_{t}^{\varepsilon}+Y_{1}<a\}} \int_{{y-\delta}-X_{t}^{\varepsilon}-Y_{1}}^{y+\delta
-X_{t}^{\varepsilon}-Y_{1}}
s_{\varepsilon}(w)\,\mathrm{d}w \biggr)
\\
&&\quad={\mathrm{e}^{-\lambda_{\varepsilon}t}}\bigl(t^{2}/2\bigr)\bbe \biggl(\int
_{(a-\underline
{X}_{t}^{\varepsilon},b-\bar{X}_{t}^{\varepsilon})^{c}} \bigl\{\varphi (v)+\|\varphi\|_{\mathrm{ Lip}}U_{t}^{\varepsilon}
\bigr\}s_{\varepsilon}(v) \int_{{y-\delta}-X_{t}^{\varepsilon}-v}^{y+\delta-X_{t}^{\varepsilon
}-v}
s_{\varepsilon}(w)\,\mathrm{d}w \,\mathrm{d}v \biggr).
\end{eqnarray*}
Next, changing variables and applying Fubini,
%
%
\begin{eqnarray} \label{Eq:UpBnd2}
A_{2,2}&=&\int_{{y-\delta}}^{y+\delta}
{{\mathrm{e}^{-\lambda_{\varepsilon
}t}}2^{-1} t^{2} \bbe \biggl(\int
_{(a-\underline{X}_{t}^{\varepsilon},b-\bar
{X}_{t}^{\varepsilon})^{c}} \bigl\{\varphi(v)+\|\varphi\|_{\mathrm{
Lip}}U_{t}^{\varepsilon}
\bigr\}s_{\varepsilon}(v)s_{\varepsilon
}\bigl(u-X_{t}^{\varepsilon}-v
\bigr)\,\mathrm{d}v \biggr)}\,\mathrm{d}u
\nonumber
\\[-4pt]
\\[-12pt]
\nonumber
&=:&\int_{{y-\delta}}^{y+\delta}\bar{P}^{2,2}_{t}(u)
\,\mathrm{d}u.
\end{eqnarray}
In order to find a lower bound for $A_{2}(t)$, note that
\begin{eqnarray*}
&&\varphi(X_{\tau})\mathbf{1}_{\{\tau\leq t, X_{t}\in I_{\delta
}(y),{N_{t}^{\varepsilon}=2}\}}\geq\varphi
\bigl(X_{\tau}^{\varepsilon
}+Y_{1}\bigr){
\mathbf{1}}_{\{\tau_{1}\leq \tau<\tau_{2}, X_{t}^{\varepsilon
}+Y_{1}+Y_{2}\in I_{\delta}(y),{N_{t}^{\varepsilon}=2}\}}
\\
&&\quad \geq\varphi\bigl({X}_{\tau}^{\varepsilon}+Y_{1}\bigr){
\mathbf{1}}_{\{
Y_{1}+\underline{X}^{\varepsilon}_{t}\geq b\ \mathrm{or}\ \bar
{X}^{\varepsilon}_{t}+Y_{1}\leq a\}}\mathbf{1}_{\{\bar{X}^{\varepsilon
}_{t}<b, \underline{X}^{\varepsilon}_{t}>a\}}\mathbf{1}_{\{
X_{t}^{\varepsilon}+Y_{1}+Y_{2}\in I_{\delta}(y),{N_{t}^{\varepsilon
}=2}\}}.
\end{eqnarray*}
Using the previous inequality and the lower bound $\varphi(X_{\tau
}^{\varepsilon}+Y_{1})\geq\varphi(Y_{1})-\|\varphi\|_{\mathrm{
Lip}}U_{t}^{\varepsilon}$ {together with the independence of
$N^{\varepsilon}$, the $Y_{i}$'s, and $X^{\varepsilon}$}, it follows that
\begin{eqnarray*}
A_{2}(t)&\geq&{\mathrm{e}^{-\lambda_{\varepsilon}t}} \frac{(\lambda_{\varepsilon}t)^{2}}{2}\bbe \bigl( \bigl\{
\varphi(Y_{1})-\|\varphi\|_{\mathrm{
Lip}}U_{t}^{\varepsilon}
\bigr\}\mathbf{1}_{\{Y_{1}\in
(a-\bar{X}^{\varepsilon}_{t},b-\underline{X}^{\varepsilon
}_{t})^{c},\bar{X}^{\varepsilon}_{t}<b, \underline{X}^{\varepsilon
}_{t}>a,X_{t}^{\varepsilon}+Y_{1}+Y_{2}\in I_{\delta}(y)\}} \bigr)
\\
&=:&\int_{{y-\delta}}^{y+\delta} \underline{P}^{2}_{t}(u)
\,\mathrm{d}u,
\end{eqnarray*}
where $\underline{P}^{2}_{t}(u)$ is defined as
\[
{{\mathrm{e}^{-\lambda_{\varepsilon}t}}2^{-1} t^{2} \bbe \biggl({
\mathbf{1}}_{\{\bar{X}^{\varepsilon}_{t}<b,
\underline{X}^{\varepsilon}_{t}>a\}}\int_{(a-\bar{X}_{t}^{\varepsilon
},b-\underline{X}_{t}^{\varepsilon})^{c}} \bigl\{ \varphi(v)-\|\varphi
\|_{\mathrm{ Lip}}U_{t}^{\varepsilon} \bigr\} s_{\varepsilon}(v)s_{\varepsilon}
\bigl(u-X_{t}^{\varepsilon}-v\bigr)\,\mathrm{d}v \biggr).}
\]
As it will be proved in Lemma \ref{PrfEstRmd0} below, {$\bar
{P}_{t}^{2,2}(u)$ and $\underline{P}_{t}^{2}(u)$ are such that}
%
%
\begin{eqnarray}
\label{EstErA} \lim_{t\to 0}\sup_{u\in(a+\delta_{0},b-\delta_{0})} \biggl
\llvert \frac
{1}{t^{2}} {\bar{P}_{t}^{2,2}(u)}-
\frac{1}{2}\int_{(a,b)^{c}} \varphi (v)s_{\varepsilon}(v)s_{\varepsilon}(u-v)\,\mathrm{d}v
\biggr\rrvert& =&0,
\\
\lim_{t\to 0}\sup_{u\in(a+\delta_{0},b-\delta_{0})} \biggl\llvert
\frac
{1}{t^{2}} {\underline{P}_{t}^{2}(u)}-
\frac{1}{2}\int_{(a,b)^{c}}\varphi(v)s_{\varepsilon}(v)s_{\varepsilon}(u-v)\,\mathrm{d}v
\biggr\rrvert &=&0. \label{EstErB}
\end{eqnarray}
Using (\ref{Eq:UpBnd1}), (\ref{EstErB}) and the corresponding bounds
for $\bar{P}^{2,1}_{t}(u)$ and $\bar{P}^{2,3}_{t}(u)$, it follows
that the nonnegative function $P_{t}^{2}(u)$ defined in (\ref{WAkIPk})
is such that
%
%
\begin{eqnarray}
\label{Eq:LmtCnvA2} {\lim_{t\to0}}\sup_{u\in(a+\delta_{0},b-\delta_{0})}
\biggl\llvert \frac
{1}{t^{2}}P_{t}^{2}(u)-
\frac{1}{2}\int_{(a,b)^{c}} \varphi(v) s_{\varepsilon}(v)s_{\varepsilon}(u-v)\,\mathrm{d}v
\biggr\rrvert =0.
\end{eqnarray}

(4) \textit{Three or more big jumps.} As before, we have the
following {bound
\begin{eqnarray*}
0&\leq&\bbe \bigl(\varphi(X_{\tau})\mathbf{1}_{\{\tau\leq t,X_{t}\in
I_{\delta}(y),N_{t}^{\varepsilon}=n\}} \bigr) \leq \|
\varphi\|_{\infty
}\bbe (\mathbf{1}_{\{X^{\varepsilon}_{t}+\sum_{i=1}^{n}Y_{i}\in
I_{\delta}(y),N_{t}^{\varepsilon}=n\}} )
\\
&=& \|\varphi\|_{\infty}\bbp\bigl(N_{t}^{\varepsilon}=n\bigr)
\int_{{y-\delta}}^{y+\delta} \bbe \bigl( s^{* n}_{\varepsilon}
\bigl(u-X_{t}^{\varepsilon}\bigr) \bigr)\,\mathrm{d}u.
\end{eqnarray*}
Using the previous inequality and (\ref{WAkIPk}), we have}
\begin{eqnarray*}
A_{3}(t)= \int_{{y-\delta}}^{y+\delta}P_{t}^{3}(u)
\,\mathrm{d}u\leq\int_{{y-\delta}}^{y+\delta} \Biggl[\sum
_{n=3}^{\infty} {\mathrm{e}^{-\lambda
_{\varepsilon}t}} \frac{t^{n}}{n!} \|
\varphi\|_{\infty} \bbe \bigl( s^{* n}_{\varepsilon}
\bigl(u-X_{t}^{\varepsilon}\bigr) \bigr) \Biggr]\,\mathrm{d}u=:\int
_{{y-\delta}}^{y+\delta} \bar{P}^{3}_{t}(u)
\,\mathrm{d}u.
\end{eqnarray*}
{Since} $\|s_{\varepsilon}^{*n}\|_{\infty}\leq\lambda
^{n-1}_{\varepsilon}a_{\varepsilon}$,
%
%
\begin{equation}
\label{Eq:DfnP3} \sup_{u}\bar{P}_{t}^{3}(u)
\leq{\mathrm{e}^{-\lambda_{\varepsilon}t}} a_{\varepsilon}\|\varphi\|_{\infty}\sum
_{n=3}^{\infty} \frac{t^{n}}{n!}\lambda_{\varepsilon}^{n-1}
\leq C(\varepsilon) t^{3}
\end{equation}
for some constant $C(\varepsilon)<\infty$, and we conclude that $0\leq
P_{t}^{3}(u)\leq C(\varepsilon) t^{3}$ for {$\calL$-}a.e. $u$.

Putting the four {previous} steps together, we conclude that
$\bbe (\varphi(X_{\tau})\mathbf{1}_{\{\tau\leq t, X_{t}\in({y-\delta
},y+\delta)\}} )=\int_{{y-\delta}}^{y+\delta} P_{t}(u)\,\mathrm{d}u$,
for a function $P_{t}(u)$ such that
\[
\lim_{t\to 0}{\sup_{u\in(a+\delta_{0},b-\delta_{0})}} \biggl\llvert
\frac
{1}{t^{2}}P_{t}(u)-\frac{1}{2}\int_{(a,b)^{c}}
\varphi(v)s_{\varepsilon
}(v)s_{\varepsilon}(u-v)\,\mathrm{d}v\biggr\rrvert =0.
\]
Finally, it is easy to see that for any $u\in(a+\delta_{0},b-\delta
_{0})$ and $a<0<b$, there exists an $\varepsilon_{0}>0$ small enough
such that
$\int_{(a,b)^{c}} \varphi(v) s_{\varepsilon}(v)s_{\varepsilon
}(u-v)\,\mathrm{d}v=\int_{(a,b)^{c}} \varphi(v) s(v)s(u-v)\,\mathrm{d}v$,
for all $0<{\varepsilon}<\varepsilon_{0}$. This concludes the proof.
\hfill$\Box$

\begin{lmma}\label{PrfEstRmd0}
Verification of \textup{(\ref{EstErA})} and \textup{(\ref{EstErB})}.
\end{lmma}
\begin{pf}
Let $0<\varepsilon<1$ and
{$M^{\varepsilon}_{t}:=X_{t}^{\varepsilon}-{\mu_{\varepsilon} }t$ be
the martingale component of $X^{\varepsilon}$}. We shall analyze the
expressions appearing
inside the absolute values on the right-hand side of equations
(\ref{EstErA}) and (\ref{EstErB}). Define the random intervals
${{\bar I}:=(a-\underline X^\varepsilon_t, b-\bar{X}_{t}^{\varepsilon
})}$, $\underline
I:=(a-\bar{X}_{t}^{\varepsilon} ,b-\underline{X}_{t}^{\varepsilon})$,
and the corresponding limiting interval $J=(a,b)$, {under the
convention $(x,y)=\varnothing$ if $y<x$}. Denote
\begin{eqnarray*}
D_{t}^{1}(u)&= &\bbe \biggl(\int_{\bar I^c}
\bigl\{\varphi(v)+\|\varphi\|_\mathrm{
Lip}U_{t}^{\varepsilon}
\bigr\}s_{\varepsilon}(v)s_{\varepsilon
}\bigl(u-X_{t}^{\varepsilon}-v
\bigr)\,\mathrm{d}v \biggr) -\int_{J^c} \varphi(v)s_{\varepsilon}(v)s_{\varepsilon}(u-v)\,\mathrm{d}v,
\\
D_{t}^{2}(u)&=& \bbe \biggl(\mathbf{ 1}_{\{\bar{X}_{t}^{\varepsilon}<b,\underline{X}_{t}^{\varepsilon}>a\}
}\int
_{
\underline I^c} \bigl\{\varphi(v)-\|\varphi\|_\mathrm{
Lip}U_{t}^{\varepsilon}
\bigr\}s_{\varepsilon}(v)s_{\varepsilon
}\bigl(u-X_{t}^{\varepsilon}-v
\bigr)\,\mathrm{d}v \biggr)
\\
&&{} -\int_{J^c} \varphi(v)s_{\varepsilon}(v)s_{\varepsilon}(u-v)\,\mathrm{d}v.
\end{eqnarray*}
Let us first analyze $D_{t}^{1}$. Clearly,
\begin{eqnarray*}
D_{t}^{1}(u) &=& {\|\varphi\|_{\mathrm{ Lip}} \bbe
\biggl(U_{t}^{\varepsilon}\int_{\bar I^c}
s_{\varepsilon
}(v)s_{\varepsilon}\bigl(u-X_{t}^{\varepsilon}-v
\bigr) \,\mathrm{d}v \biggr)}
\\
&&{} + {\bbe \biggl(\int_{\bar I^c \setminus J^c} \varphi (v)s_{\varepsilon}(v)s_{\varepsilon}
\bigl(u-X_{t}^{\varepsilon}-v\bigr) \,\mathrm{d}v \biggr)}
\\
&&{} +\bbe \biggl(\int_{J^c} \varphi(v)s_{\varepsilon}(v)
\bigl[s_{\varepsilon}\bigl(u-X_{t}^{\varepsilon}-v\bigr)
-s_{\varepsilon}(u-v)\bigr]\,\mathrm{d}v \biggr),
\end{eqnarray*}
and, therefore, {using that $\bar{I}^{c}\setminus J^{c}\subset
(a,a-\underline{X}_{t}^{\varepsilon})\cup (b-\bar{X}_{t}^{\varepsilon
},b)$, under the convention $(-\infty,-\infty)=(\infty,\infty)=\varnothing$,}
\begin{eqnarray*}
\bigl|D_{t}^{1}(u)\bigr|&\leq& a_{\varepsilon}\lambda_{\varepsilon}
\|\varphi\| _\mathrm{ Lip} \bbe {U}_{t}^{\varepsilon} +
a_\varepsilon^2 \|\varphi\|_\infty\bbe\bigl(\bar
X^\varepsilon_t - \underline X^\varepsilon_t
\bigr) + \lambda_{\varepsilon}\|\varphi \|_{\infty}\bigl\|s'_{\varepsilon}
\bigr\|_{\infty} \bbe\bigl|X_{t}^{\varepsilon}\bigr|
\\
&\leq&\bigl(a_{\varepsilon}\lambda_{\varepsilon}\|\varphi\|_\mathrm{
Lip}+2\|
\varphi\|_{\infty}a_{\varepsilon}^{2}\bigr) \Bigl(\bbe \sup
_{s\leq t}\bigl |{M}_{s}^{\varepsilon}\bigr|+|{
\mu_{\varepsilon}}| t \Bigr) +\|\varphi\|_{\infty}\lambda_{\varepsilon}
\bigl\|s'_{\varepsilon}\bigr\|_{\infty
} \bigl(\bbe\bigl|M_{t}^{\varepsilon}\bigr|+|{
\mu_{\varepsilon}}|t \bigr).
\end{eqnarray*}
{Using the trivial inequalities $(\bbe\sup_{s\leq t}
|{M}_{s}^{\varepsilon}|)^{2}\leq\bbe\sup_{s\leq t}
({M}_{s}^{\varepsilon})^{2}$ and $(\bbe|{M}_{s}^{\varepsilon
}|)^{2}\leq\bbe({M}_{s}^{\varepsilon})^{2}$, together with} Doob's
inequality, we then get the bound
%
%
\begin{eqnarray}\label{Eq:BndKErr}
\bigl\llvert D_{t}^{1}(u)\bigr\rrvert &\leq&
\bigl[{2}a_{\varepsilon}\lambda _{\varepsilon}\|\varphi\|_\mathrm{ Lip}+4\|
\varphi\|_{\infty}a_{\varepsilon
}^{2}+\|\varphi\|_{\infty}
\lambda_{\varepsilon}\bigl\|s'_{\varepsilon}\bigr\| _{\infty} \bigr]
\sigma_{\varepsilon}t^{1/2}
\nonumber
\\[-8pt]
\\[-8pt]
\nonumber
&&{} + \bigl[a_{\varepsilon}\lambda_{\varepsilon}\|\varphi\|_\mathrm{
Lip}+2\|
\varphi\|_{\infty}a_{\varepsilon}^{2}+\|\varphi
\|_{\infty
}\lambda_{\varepsilon}\bigl\|s'_{\varepsilon}
\bigr\|_{\infty} \bigr]|{ \mu_{\varepsilon}}|t,
\end{eqnarray}
where $\sigma_{\varepsilon}^{2}:=\sigma^{2}+\int\bar{c}_{\varepsilon
}(x) x^{2} \nu(\mathrm{d}x)$. For $D_{t}^{2}(u)$, note that
\begin{eqnarray*}
D_{t}^{2}(u)&=& - \bbe \biggl(\mathbf{1}_{\{\bar{X}_{t}^{\varepsilon}\geq b\ \mathrm{or}\ \underline{X}_{t}^{\varepsilon}\leq
a\}}\int
_{\underline I^c} \varphi(v)s_{\varepsilon}(v)s_{\varepsilon}
\bigl(u-X_{t}^{\varepsilon}-v\bigr)\,\mathrm{d}v \biggr)
\\
&&{} + \|\varphi\|_\mathrm{
Lip}\bbe \biggl(\mathbf{ 1}_{\{\bar{X}_{t}^{\varepsilon}<b,\underline{X}_{t}^{\varepsilon}>a\}}
U_{t}^{\varepsilon}\int_{
\underline I^c}
s_{\varepsilon}(v)s_{\varepsilon}\bigl(u-X_{t}^{\varepsilon}-v
\bigr)\,\mathrm{d}v \biggr)
\\
&&{} + \bbe \biggl(\int_{\underline I^c} \varphi(v)s_{\varepsilon
}(v)s_{\varepsilon}
\bigl(u-X_{t}^{\varepsilon}-v\bigr)\,\mathrm{d}v \biggr) -\int
_{J^c} \varphi(v)s_{\varepsilon}(v)s_{\varepsilon}(u-v)\,\mathrm{d}v.
\end{eqnarray*}
Defining $c=|a|\wedge b$ and following the same steps as above, it is
easy to verify that $|D_{t}^{2}(u)|$ admits the following upper bound:
%
%
\begin{eqnarray}\label{Eq:BndKErr2}
{\bigl\llvert D_{t}^{2}(u)\bigr\rrvert }&\leq&\|\varphi
\|_{\infty}\lambda _{\varepsilon}a_{\varepsilon}\bbp\bigl(U^{\varepsilon}_{t}
\geq c\bigr)+a_{\varepsilon}\lambda_{\varepsilon}\|\varphi\|_\mathrm{ Lip}
\bbe {U}_{t}^{\varepsilon}+2a_{\varepsilon}^{2} \|\varphi
\|_{\infty}\bbe U_{t}^{\varepsilon} +\lambda_{\varepsilon}\|
\varphi\|_{\infty}\bigl\|s'_{\varepsilon}\bigr\|_{\infty
}
\bbe\bigl|X_{t}^{\varepsilon}\bigr|
\nonumber
\\
& \leq &\|\varphi\|_{\infty}\lambda_{\varepsilon}a_{\varepsilon
}C_{1}(c,
\varepsilon) t+ \bigl[ 2a_{\varepsilon}\lambda_{\varepsilon}\| \varphi
\|_\mathrm{ Lip}+4\|\varphi\|_{\infty}a_{\varepsilon}^{2}+\|
\varphi \|_{\infty}\lambda_{\varepsilon}\bigl\|s'_{\varepsilon}
\bigr\|_{\infty} \bigr]\sigma _{\varepsilon}t^{1/2}
\\
&&{} + \bigl[a_{\varepsilon}\lambda_{\varepsilon}\|\varphi\|_\mathrm{
Lip}+{2}
\|\varphi\|_{\infty}a_{\varepsilon}^{2}+\|\varphi
\|_{\infty
}\lambda_{\varepsilon}\bigl\|s'_{\varepsilon}
\bigr\|_{\infty} \bigr]|{ \mu_{\varepsilon}}|t,\nonumber
\end{eqnarray}
where we had used the tail probability bound in (\ref{KyInBn}).
\end{pf}

\subsection{Proof of Proposition \texorpdfstring{\protect\ref{Th:y_out}}{2.1}}

We use the notation introduced at the beginning of Section~\ref{SbSec:PrMR} above and, as before, we assume without loss of generality
that $\varphi$ is nonnegative. As it was done in (\ref{T1-T4}), by
{partitioning the space into the different values that
$N_{t}^{\varepsilon}$ can take on,} we can decompose {$\bbe(\varphi
(X_{\tau})\mathbf{
1}_{\{X_{t}\in I_{\delta}(y)\}})$} into three terms: no big
jumps, one big jump, and two or more big jumps. These terms can in
turn be expressed as integrals of the form (\ref{WAkIPk}) using a
procedure similar to (\ref{WAkIPkb}).
The term with no big jumps is such that
\begin{eqnarray*}
\int_{{y-\delta}}^{y+\delta} P_{t}^{0}(u)\,\mathrm{d}u&:=&
\bbe\bigl( \varphi(X_{\tau
})\mathbf{1}_{\{X_{t}\in
I_{\delta}(y),N_{t}^{\varepsilon}=0\}}\bigr)
\\
& \leq& \|\varphi\|_\infty\bbp \bigl(X^\varepsilon_t
\in I_\delta (y),N_{t}^{\varepsilon}=0 \bigr)
\\
& =& \int_{{y-\delta}}^{y+\delta} \mathrm{e}^{-\lambda_{\varepsilon}t}\|\varphi\|
_\infty f^\varepsilon_t(u)\,\mathrm{d}u,
\end{eqnarray*}
{which yields an upper bound for $P_{t}^{0}(u)$ of the form $\bar
{P}^{0}_{t}(u):=\mathrm{e}^{-\lambda_{\varepsilon}t}\|\varphi\|_\infty
f^\varepsilon_t(u)$. Using\vspace*{1pt} (\ref{KyInBn})(ii), we can further upper
bound $\bar{P}^{0}_{t}(u)$ by $\|\varphi\|_{\infty} c_{2}(c, \varepsilon
) t^{2}$ uniformly in $(a-\delta_{0},b+\delta_{0})^{c}$}.
The term with two or more big jumps can be bounded similarly to the
term with three or more big jumps in the previous section. {Concretely,
this term satisfies
\begin{eqnarray*}
\int_{{y-\delta}}^{y+\delta}P_{t}^{2}(u)
\,\mathrm{d}u&:=& \bbe \bigl( \varphi(X_{\tau})\mathbf{ 1}_{\{X_{t}\in I_{\delta}(y), N_{t}^{\varepsilon}\geq 2\}} \bigr)
\\
&\leq& \int_{{y-\delta}}^{y+\delta} \Biggl[\sum
_{n=2}^{\infty} {\mathrm{e}^{-\lambda
_{\varepsilon}t}} \frac{t^{n}}{n!} \|
\varphi\|_{\infty} \bbe \bigl( s^{* n}_{\varepsilon}
\bigl(u-X_{t}^{\varepsilon}\bigr) \bigr) \Biggr]\,\mathrm{d}u=:\int
_{{y-\delta}}^{y+\delta} \bar{P}^{2}_{t}(u)
\,\mathrm{d}u,
\end{eqnarray*}
and, using that $\|s_{\varepsilon}^{*n}\|_{\infty}\leq\lambda
^{n-1}_{\varepsilon}a_{\varepsilon}$, we can further upper bound $ \bar
{P}^{2}_{t}(u)$ by $C(\varepsilon)t^{2}$ for a constant $C(\varepsilon
)<\infty$.}
The term with exactly one jump is decomposed as follows:
\begin{eqnarray*}
\int_{{y-\delta}}^{y+\delta}P_{t}^{1}(u)
\,\mathrm{d}u&:=&\bbe \bigl( \varphi (X_{\tau})\mathbf{1}_{\{X_{t}\in
I_{\delta}(y), N_{t}^{\varepsilon}=1\}} \bigr)
\\
&=& \bbe \bigl( \varphi\bigl(X^\varepsilon_{\tau}\bigr){
\mathbf{1}}_{\{X_{t}\in I_{\delta
}(y); \tau<
\tau_1; N_{t}^{\varepsilon}=1\}} \bigr)+ \bbe \bigl(\varphi(X_{\tau}){
\mathbf{1}}_{\{X_{t}\in I_{\delta}(y);
\tau\geq\tau_1;N_{t}^{\varepsilon}=1\}} \bigr),
\end{eqnarray*}
where $\tau_1$ is the time of the first big jump.
Out of these two terms, the first one {satisfies
\begin{eqnarray*}
\bbe \bigl( \varphi\bigl(X^\varepsilon_{\tau}\bigr){
\mathbf{1}}_{\{X_{t}\in I_{\delta
}(y); \tau<
\tau_1; N_{t}^{\varepsilon}=1\}} \bigr)&\leq&\|\varphi\|_\infty \bbp\bigl[\exists s
\in[0,t]\dvt X^\varepsilon_s \notin D; X^\varepsilon_t
+ Y_1 \in I_\delta; N_{t}^{\varepsilon}=1\bigr]
\\
&=& \mathrm{e}^{- \lambda_{\varepsilon}t}\lambda_{\varepsilon}t\|\varphi\|_\infty \bbp\bigl[
\exists s\in[0,t]\dvt X^\varepsilon_s \notin D;
X^\varepsilon_t + Y_1 \in I_\delta\bigr]
\\
&\leq&\int_{{y-\delta}}^{y+\delta}\mathrm{e}^{- \lambda_{\varepsilon}t}t\| \varphi
\|_\infty \bbe \bigl[{\mathbf{1}}_{U^\varepsilon_t \geq c} s_\varepsilon
\bigl(u-X_{t}^{\varepsilon}\bigr) \bigr]\,\mathrm{d}u,
\end{eqnarray*}
where the integrand $\bar{P}^{1,1}_{t}(u):=\mathrm{e}^{- \lambda_{\varepsilon
}t}t\|\varphi\|_\infty \bbe [{\mathbf{1}}_{U^\varepsilon_t \geq c}
s_\varepsilon(u-X_{t}^{\varepsilon})  ]$ is uniformly bounded by
$\|\varphi\|_\infty
a_\varepsilon C_1(c,\varepsilon)t^{2}$.
As} for the second term
\[
\bbe \bigl( \varphi(X_{\tau})\mathbf{1}_{\{X_{t}\in I_{\delta}(y);
\tau\geq\tau_1; N_{t}^{\varepsilon}=1\}} \bigr) = \bbe
\bigl( \varphi\bigl(X^\varepsilon_{\tau} + Y_1\bigr){
\mathbf{1}}_{\{
X_{t}^{\varepsilon}+Y_{1}\in
I_{\delta}(y); \tau\geq\tau_1;N_{t}^{\varepsilon}=1\}} \bigr),
\]
it can be bounded from above {by
\begin{eqnarray*}
&&\bbe\bigl( \varphi(Y_1)\mathbf{1}_{\{X^\varepsilon_t + Y_1\in
I_{\delta}(y);N_{t}^{\varepsilon}=1\}}\bigr) +\|\varphi
\|_{\mathrm{Lip}} \bbe \bigl( U^\varepsilon_t {
\mathbf{1}}_{\{X^\varepsilon_t + Y_1\in
I_{\delta}(y);N_{t}^{\varepsilon}=1\}}\bigr)
\\
&&\quad= \int_{{y-\delta}}^{y+\delta} \bigl\{\mathrm{e}^{- \lambda_{\varepsilon}t}t \bbe
\bigl[ \varphi\bigl(u-X_{t}^{\varepsilon}\bigr) s_\varepsilon
\bigl(u-X_{t}^{\varepsilon
}\bigr) \bigr] +\mathrm{e}^{- \lambda_{\varepsilon}t}t \|\varphi
\|_{\mathrm{Lip}} \bbe \bigl[U^\varepsilon_t s_\varepsilon
\bigl(u-X_{t}^{\varepsilon}\bigr) \bigr] \bigr\} \,\mathrm{d}u
\\
&&\quad\leq\int_{{y-\delta}}^{y+\delta} \bigl\{t\varphi(u)
s_\varepsilon(u) +t\bigl( \|\varphi\|_{\mathrm{Lip}} a_\varepsilon+ \|
\varphi\|_\infty\bigl\|s'_\varepsilon\bigr\|_\infty
\bigr) \bbe\bigl[\bigl|X^\varepsilon_t\bigr|\bigr] + t \|\varphi
\|_{\mathrm{Lip}} a_\varepsilon \bbe\bigl[U^\varepsilon_t
\bigr] \bigr\} \,\mathrm{d}u.
\end{eqnarray*}
Similarly, this can be bounded from below by
\begin{eqnarray*}
&&\hspace*{-4pt} \bbe\bigl(\varphi(Y_1)\mathbf{1}_{\{X^\varepsilon_{t} + Y_1\in I_{\delta}(y);
\bar X^\varepsilon_t < b; \underline X^\varepsilon_t > a;
N_{t}^{\varepsilon}=1\}}\bigr) -\|\varphi
\|_{\mathrm{Lip}} \bbe\bigl( U^\varepsilon_t {
\mathbf{1}}_{\{X^\varepsilon_t + Y_1\in
I_{\delta}(y); N_{t}^{\varepsilon}=1\}}\bigr)
\\
&&\hspace*{-4pt}\quad=\int_{y-\delta}^{y+\delta} \bigl\{\mathrm{e}^{- \lambda_{\varepsilon}t}t \bbe
\bigl(\varphi\bigl(u-X_{t}^{\varepsilon}\bigr)s_{\varepsilon}
\bigl(u-X_{t}^{\varepsilon
}\bigr)\mathbf{1}_{\{
\bar X^\varepsilon_t < b, \underline X^\varepsilon_t > a\}}\bigr)
-\mathrm{e}^{-
\lambda_{\varepsilon}t}t \|\varphi\|_{\mathrm{Lip}} \bbe\bigl( U^\varepsilon_t
s_{\varepsilon}\bigl(u-X_{t}^{\varepsilon}\bigr)\bigr) \bigr\} \,\mathrm{d}u
\\
&&\hspace*{-4pt}\quad \geq\int_{y-\delta}^{y+\delta} \bigl\{ t\varphi(u)
s_\varepsilon(u) -\|\varphi\|_\infty a_\varepsilon\lambda
_{\varepsilon}t^{2}-t\bigl( \|\varphi\|_{\mathrm{Lip}}
a_\varepsilon+ \|\varphi\|_\infty\bigl\|s'_\varepsilon
\bigr\|_\infty\bigr) \bbe\bigl[\bigl|X^\varepsilon_t\bigr|\bigr]
\\
&&\hspace*{24pt}\qquad{} -t \|\varphi\|_{\mathrm{Lip}} a_\varepsilon \bbe\bigl[U^\varepsilon_t
\bigr] - t\|\varphi\|_\infty a_\varepsilon C_1(c,
\varepsilon) t \bigr\}\,\mathrm{d}u.
\end{eqnarray*}
To} conclude, we estimate $\bbe[|X^\varepsilon_t|]$
and $\bbe[U^\varepsilon_t]$ as in the proof of Lemma \ref{PrfEstRmd0}
above.

\section{Proofs of Section \texorpdfstring{\protect\ref{estgen}}{3}}\label{Ap:PrfErrBnd}
In this part, we provide the building blocks to develop an upper bound
for the remainder $\mathcal{R}_{t}(u)$ appearing in (\ref{Eq:FMR}).

\subsection{Proof of Lemma \texorpdfstring{\protect\ref{Lmm:KTPBnd}}{3.1}}

{Let us first assume that ${\mu}\geq0$ so that
$X_{t}:=M_{t}+{\mu}t$ is a submartingale.} By Doob's inequality, for
all $c>0$,
%
%
\begin{equation}
\mathbb P \Bigl(\sup_{s\leq t} {X_s} \geq\eta
\Bigr) = \mathbb P \Bigl(\sup_{s\leq t} \mathrm{e}^{c {X_s}} \geq
\mathrm{e}^{c\eta} \Bigr) \leq\frac{\mathbb E[\mathrm{e}^{c {X_t}}]}{\mathrm{e}^{c\eta}} = \mathrm{e}^{t \psi(c) - c \eta
}\label{doob}
\end{equation}
with $\psi(c) = {{\mu}c}+\frac{\sigma^2 c^2}{2} + \int_{|z|\leq
\varepsilon} (\mathrm{e}^{cz} - 1
- cz)\nu(\mathrm{d}z)$. Minimizing\vspace*{1pt} the right-hand side over all $c>0$, we get,
as in
R\"uschendorf and
Woerner \cite{ruschendorf02} {(see page~87 therein)},
%
%
\begin{equation}
\inf_{c>0} \mathrm{e}^{t \psi(c) - c \eta} = \exp \biggl(-t \int
_{{\psi'(0)}}^{\eta
/t} \tau(z)\,\mathrm{d}z \biggr)=\exp \biggl(-t \int
_{{{\mu}}}^{\eta/t} \tau (z)\,\mathrm{d}z \biggr),\label{doob2}
\end{equation}
where {we are taking $t<\eta/{\mu}$} and ${\tau\dvtx  [0,\infty)\to\bbr}$ is
the inverse function
of
%
%
\begin{equation}
\label{Eq:Psiprime} \psi'(x) = {\mu}+\sigma^2 x + \int
_{|z|\leq\varepsilon} z\bigl(\mathrm{e}^{zx}-1\bigr)\nu(\mathrm{d}z).
\end{equation}
As in Houdr\'e \cite{Houdre:2002}, note that, for $x \geq0$,
\begin{eqnarray*}
0\leq\int_{|z|\leq\varepsilon} z\bigl(\mathrm{e}^{zx}-1\bigr)\nu(\mathrm{d}z)&\leq&
\int_{|z|\leq
\varepsilon} |z|\bigl(\mathrm{e}^{|z|x}-1\bigr)\nu(\mathrm{d}z)\leq\int
_{|z|\leq\varepsilon} |z|\sum_{k=1}^{\infty}
\frac{(|z|x)^{k}}{k!}\nu(\mathrm{d}z)
\\
&\leq&\int_{|z|\leq\varepsilon} |z|^{2}\nu(\mathrm{d}z) \sum
_{k=1}^{\infty} \frac{\varepsilon^{k-1} x^{k}}{k!}=\int
_{|z|\leq\varepsilon} |z|^{2}\nu (\mathrm{d}z) \frac{1}{\varepsilon}
\bigl(\mathrm{e}^{\varepsilon x}-1 \bigr).
\end{eqnarray*}
From the previous inequality, for $x\geq 0$,
\[
0\leq\psi'(x) \leq{\mu}+\sigma^2 x +\int
_{|z|\leq\varepsilon} |z|^{2}\nu(\mathrm{d}z) \frac{1}{\varepsilon}
\bigl(\mathrm{e}^{\varepsilon x}-1 \bigr) \leq {\mu}+\frac{\mathrm{e}^{\varepsilon x}-1}{\varepsilon}
\sigma_\varepsilon^2,
\]
where we used the fact that $\sigma^{2}_{\varepsilon}=\sigma^{2}+\int_{|z|\leq\varepsilon} |z|^{2}\nu(\mathrm{d}z)$.
This implies that
\[
\tau(z)\geq\frac{1}{\varepsilon}\log \biggl\{1 + \frac{{z-{\mu}}}{\sigma_\varepsilon^2} \varepsilon
\biggr\},
\]
and therefore, substituting this into \eqref{doob} and \eqref{doob2}
and using that $v\ln(v)\leq(1+v)\ln(1+v)$ and $\mathrm{e}^{-v\log v}\leq
\mathrm{e}^{\mathrm{e}^{-1}}$ for all $v>0$, we have
\begin{eqnarray*}
\bbp \Bigl[\sup_{s\leq t} {X_s} \geq \eta \Bigr] &
\leq &\exp \biggl\{-\frac{t\sigma_\varepsilon^2}{\varepsilon^2}\int_{0}^{{\varepsilon{(\eta-{\mu}t)}}/{(t\sigma_\varepsilon^2)}}
\log(1+s)\,\mathrm{d}s \biggr\}
\\
&=&\exp \biggl\{-\frac{t\sigma_\varepsilon^2}{\varepsilon^2} \biggl( \biggl(1+\frac{\varepsilon
{(\eta-{\mu}t)}}{t\sigma_\varepsilon^2} \biggr)
\log \biggl(1+\frac
{\varepsilon
{(\eta-{\mu}t)}}{t\sigma_\varepsilon^2} \biggr)-\frac{\varepsilon
{(\eta-{\mu}t)}}{t\sigma_\varepsilon^2} \biggr) \biggr\}
\\
&\leq &\exp \biggl\{-\frac{{\eta-{\mu}t}}{\varepsilon}\log \biggl(\frac{
\varepsilon{(\eta-{\mu}t)}}{\mathrm{e} \sigma_\varepsilon^2 t} \biggr)
\biggr\}
\\
&\leq& t^{{\eta}/{\varepsilon}}\mathrm{e}^{({{\mu}}/{\varepsilon
})\mathrm{e}^{-1}}\exp \biggl\{-\frac{{\eta-{\mu}t}}{\varepsilon}\log
\biggl(\frac{
\varepsilon{(\eta-{\mu}t)}}{\mathrm{e} \sigma_\varepsilon^2} \biggr) \biggr\},
\end{eqnarray*}
The above inequality proves the statement (2)(i) for the case ${\mu
}=0$. Next, it is easy to check that the function $u\to(u/\varepsilon
)\log(\varepsilon u/\mathrm{e}\sigma_{\varepsilon}^{2})$ is strictly convex in
$(0,\infty)$ and reaches its global minimum value of $-\sigma
_{\varepsilon}^{2}/\varepsilon^{2}$ at $u=\sigma_{\varepsilon
}^{2}/\varepsilon$. Hence, whenever $\eta- {\mu}t\geq 0$,
\[
\bbp \Bigl[\sup_{s\leq t} {X_s} \geq \eta \Bigr]
\leq t^{{\eta}/{\varepsilon}}\mathrm{e}^{({{\mu}}/{\varepsilon})\mathrm{e}^{-1}}\mathrm{e}^{{\sigma_\varepsilon^2}/{\varepsilon^{2}}},
\]
which proves the statement (1) for ${\mu}\geq 0$. Also, if ${\mu}>0$,
$t<\eta/{\mu}$, and $\eta<\sigma_{\varepsilon}^{2}/\varepsilon$, we
have that
\[
\exp \biggl\{-\frac{{\eta-{\mu}t}}{\varepsilon}\log \biggl(\frac{
\varepsilon{(\eta-{\mu}t)}}{\mathrm{e} \sigma_\varepsilon^2} \biggr) \biggr\}
\leq \sup_{0\leq u\leq \eta}\exp \biggl\{-\frac{u}{\varepsilon}\log \biggl(
\frac{
\varepsilon{u}}{\mathrm{e} \sigma_\varepsilon^2} \biggr) \biggr\}= \biggl(\frac
{\mathrm{e}\sigma_\varepsilon^2}{\varepsilon\eta}
\biggr)^{{\eta
}/{\varepsilon}},
\]
which proves the statement (2)(ii).
Finally, we consider the case ${\mu}<0$. In that case, obviously,
$M_{t}+{\mu}t \leq M_{t}$ and
\[
\mathbb P \Bigl(\sup_{s\leq t} (M_s+ {\mu} s ) \geq
\eta \Bigr) \leq\mathbb P \Bigl(\sup_{s\leq t} M_s
\geq\eta \Bigr)\leq t^{{\eta}/{\varepsilon}} \biggl(\frac{e\sigma_\varepsilon
^2}{\varepsilon\eta}
\biggr)^{{\eta}/{\varepsilon}}\leq t^{{\eta}/{\varepsilon}} \mathrm{e}^{{\sigma_\varepsilon^2}/{\varepsilon^{2}}} ,
\]
where in the second inequality we used the case (2)(i) with ${\mu}=0$
that was proved above. The previous inequality proves the bounds (2)(i)
and (1) for ${\mu}<0$.

\subsection{Proof of Theorem \texorpdfstring{\protect\ref{Th:KRFFEB}}{3.1}}

To prove the estimate (\ref{TenOvBnd}) for the remainder $\mathcal
{R}_{t}(y)$, we analyze each of the four terms in (\ref{T1-T4})
contributing to it.

(\textit{No big jump}). The first component of the error is due to
{${P}_{t}^{0}$ which, as seen in (\ref{Eq:FrsTrmRmd}), can be bounded by
\[
\mathrm{e}^{(0)}(0,y,t):=\mathrm{e}^{-\lambda_{\varepsilon} t}\|\varphi\|_\infty \sup
_{0< u \leq t} \sup_{x\in(a,b)^c} f_{u}^{\varepsilon
}(y-x)=\mathrm{e}^{-\lambda_{\varepsilon} t}
\|\varphi\|_\infty \sup_{0< u \leq t} \sup
_{z\in(y-b,y - a)^c}f_{u}^{\varepsilon}(z).
\]
Next, recalling the notation $\Delta_y= (b-y)\wedge(y - a)>0$ and
employing our hypothesis that $X_{t}^{\varepsilon}$ has unimodal
distribution, we can further apply the bound (\ref{Eq:Bnd2}) to get}
\begin{eqnarray*}
\mathrm{e}^{(0)}(0,y,t)\leq \mathrm{e}^{-\lambda_{\varepsilon}
t}\frac{4\|\varphi\|_\infty}{\Delta_y}\sup
_{0<u\leq t}\bbp \biggl[\bigl|X_{u}^{\varepsilon}\bigr|\geq
\frac{\Delta_y}{2} \biggr]\leq \frac{8 \mathrm{e}^{-\lambda_{\varepsilon} t}\|\varphi\|_\infty}{\Delta_y} {C}(\Delta_y/4,
\varepsilon)t^{{\Delta_y}/{(4\varepsilon)}}
\end{eqnarray*}
for $t<t_{0}(\varepsilon,\Delta_y/2) \wedge t_{1}(\varepsilon,\Delta_y/2)$.

(\textit{One big jump}).
There are two sub-components to the error {in} this case. The first is
due to $\bar{P}^{1,1}$ in (\ref{Eq:DefnP11}). This term can be bounded by
\begin{eqnarray*}
\mathrm{e}^{(1,1)}(0,y,t)&:=&\|\varphi\|_\infty \mathrm{e}^{-\lambda_{\varepsilon}
t}t\bbe
\bigl(\mathbf{ 1}_{\{{U}_{t}^{\varepsilon}\geq c\}}s_{\varepsilon
}\bigl(y-X_{t}^{\varepsilon}
\bigr) \bigr)
\\
&\leq& \|\varphi\|_\infty \mathrm{e}^{-\lambda_{\varepsilon} t} ta_{\varepsilon} \bbp
\bigl({U}_{t}^{\varepsilon}\geq c\bigr)
\\
&\leq& 2\|\varphi\|_\infty \mathrm{e}^{-\lambda_{\varepsilon} t} a_{\varepsilon} {{C}(c/2,
\varepsilon) t^{1+{c}/{(2\varepsilon)}}}
\end{eqnarray*}
for $t<t_{0}(\varepsilon,c/2)$.
The other sub-component is due to $\bar{P}^{1,2}$ in (\ref
{Eq:DefnP12}), which can be bounded, {for $t<t_{0}(\varepsilon,\Delta
_y/2)$}, as follows:
\begin{eqnarray*}
\mathrm{e}^{(1,2)}(0,y,t)&:=&\|\varphi\|_\infty \mathrm{e}^{-\lambda_{\varepsilon}
t}t\bbe
\bigl(\mathbf{ 1}_{\{\bar{X}_{t}^{\varepsilon}-X_{t}^{\varepsilon}+y\geq b \ \mathrm{or}\  \underline{X}_{t}^{\varepsilon}-X_{t}^{\varepsilon}+y\leq a
\}}s_{\varepsilon}\bigl(y-X_{t}^{\varepsilon}
\bigr) \bigr)
\\
& \leq&\|\varphi\|_\infty \mathrm{e}^{-\lambda_{\varepsilon} t}t a_{\varepsilon} \bbp
\biggl(\sup_{u\leq t}\bigl|{X}_{u}^{\varepsilon}\bigr|\geq
\frac{\Delta
_y}{2} \biggr)\leq 2\|\varphi\|_\infty \mathrm{e}^{-\lambda_{\varepsilon} t}
a_{\varepsilon} C(\Delta_y/4,\varepsilon) t^{1+{{\Delta_y}/{(4\varepsilon)}}}.
\end{eqnarray*}

(\textit{Three or more big jumps}).
This component {can be} bounded as in (\ref{Eq:DfnP3}):
%
%
\begin{eqnarray}\label{Eq:DfnP3B}
\mathrm{e}^{(3)}(0,y,t)&:=&\|\varphi\|_\infty \mathrm{e}^{-\lambda_{\varepsilon} t}
\sum_{n=3}^{\infty} \frac{t^{n}}{n!} \bbe
\bigl( s^{* n}_{\varepsilon}\bigl(u-X_{t}^{\varepsilon}
\bigr) \bigr)
\nonumber
\\[-8pt]
\\[-8pt]
\nonumber
&\leq& \|\varphi\|_\infty a_{\varepsilon}\lambda_{\varepsilon}^{-1}
\bigl(1 -\mathrm{e}^{-\lambda_{\varepsilon}t}\bigl[1+\lambda_{\varepsilon}t+ (\lambda_{\varepsilon}t)^{2}/2
\bigr] \bigr).
\end{eqnarray}

(\textit{Two big jumps}).
There are three sub-components to the error in this case. From (\ref
{Eq:DfnP21}),
%
%
\begin{eqnarray}\label{Eq:DfnP23B}
\mathrm{e}^{(2,1)}(0,y,t)&:=&\|\varphi\|_\infty \mathrm{e}^{-\lambda_{\varepsilon} t}
\frac
{t^{2}}{2} \int_{-\infty}^{\infty}
s_{\varepsilon}(v) \bbe \bigl\{ \mathbf{1}_{\{{U}_{t}^{\varepsilon}\geq c\}}s_{\varepsilon
}
\bigl(y-X_{t}^{\varepsilon}-v\bigr) \bigr\} \,\mathrm{d}v
\nonumber
\\[-8pt]
\\[-8pt]
\nonumber
&\leq&\|\varphi\|_\infty \mathrm{e}^{-\lambda_{\varepsilon} t}{a_{\varepsilon
}
\lambda_{\varepsilon}} {{C}(c/2,\varepsilon)t^{2+{c}/{(2\varepsilon)}}}
\end{eqnarray}
for $t<t_{0}(\varepsilon,c)$. Similarly, from (\ref{Eq:DfnP23}),
%
%
\begin{eqnarray}\label{Eq:DfnP23B}
&&\mathrm{e}^{(2,3)}(0,y,t)\nonumber\\
&&\quad:=\|\varphi\|_\infty \mathrm{e}^{-\lambda_{\varepsilon} t}
\frac
{t^{2}}{2} \int_{-\infty}^{\infty}
s_{\varepsilon}(v) \bbe \bigl\{ \mathbf{ 1}_{\{\bar{X}_{t}^{\varepsilon}-{X}_{t}^{\varepsilon}+y\geq b \ \mathrm{or}\
\underline{X}_{t}^{\varepsilon}-{X}_{t}^{\varepsilon}+y\leq a\}
}s_{\varepsilon}
\bigl(y-X_{t}^{\varepsilon}-v\bigr) \bigr\} \,\mathrm{d}v
\\
&&\quad\leq \|\varphi\|_\infty \mathrm{e}^{-\lambda_{\varepsilon}
t}\frac{t^{2}}{2}
a_{\varepsilon}\lambda_{\varepsilon} \bbp \biggl(\sup_{u\leq t}\bigl|{X}_{u}^{\varepsilon}\bigr|
\geq \frac{\Delta
_y}{2} \biggr)\leq \|\varphi\|_\infty
\mathrm{e}^{-\lambda_{\varepsilon}
t}a_{\varepsilon}\lambda_{\varepsilon} C(\Delta_y/4,
\varepsilon) t^{2+{{\Delta_y}/{(4\varepsilon)}}}\nonumber
\end{eqnarray}
for $t<t_{0}(\varepsilon,\Delta_y/2)$.
Next, we consider {the} error due to the limits
(\ref{EstErA})--(\ref{EstErB}). These were bounded in Lemma
\ref{PrfEstRmd0}. Hence, by taking the maximum of (\ref{Eq:BndKErr})
and \eqref{Eq:BndKErr2}, after some simplification, {we get the
following expression for the error term $\mathrm{e}^{(2,2)}(0,y,t)$}:
\begin{eqnarray*}
&&\mathrm{e}^{-\lambda_{\varepsilon}t}t^{2} \biggl(\|\varphi\|_{\infty}
\lambda_{\varepsilon}a_{\varepsilon
}C(c/2,\varepsilon) t^{{c}/{(2\varepsilon)}}\\
&&\hspace*{34pt}{}+
\bigl[a_{\varepsilon}\lambda_{\varepsilon}\| \varphi\|_\mathrm{ Lip}+ 2\|
\varphi\|_{\infty}a_{\varepsilon}^{2}+\|\varphi\|_{\infty}
\lambda _{\varepsilon}\bigl\|s'_{\varepsilon}\bigr\|_{\infty} \bigr]
\biggl(\sigma_{\varepsilon}t^{1/2} + \frac{|{\mu_\varepsilon}|}{2} t\biggr)
\biggr).
\end{eqnarray*}

Finally, we also need to take into account the error due to approximating
\[
\mathrm{e}^{-\lambda_{\varepsilon}t}\frac{t^{2}}{2}\int_{(a,b)^c} \varphi(v)
s_{\varepsilon}(v)s_{\varepsilon}(y-v)\,\mathrm{d}v
\]
by $\frac{t^{2}}{2}\int_{(a,b)^c} \varphi(v) s_{\varepsilon
}(v)s_{\varepsilon}(y-v)\,\mathrm{d}v$, which is of order $\|\varphi\|_\infty
\lambda^{2}_{\varepsilon} a_{\varepsilon} t^{3}/2$. Putting all the
previous bounds together, we obtain the overall bound (\ref{TenOvBnd}).

\section{Finding the estimate $e_{f}(0,y,t)$ for the Cauchy process}
\label{cauchy.sec}

In this paragraph, our aim is to find an explicit bound for
the Cauchy process with L\'evy density $\nu(x) = \frac{c}{|x|^2}$ (and
no drift), which is used in the numerical illustrations. For
simplicity, we shall only consider the one-sided case ($a=-\infty$).
Setting $c_\varepsilon(x) = \mathbf{1}_{|x|>\varepsilon}$, we get
${\mu_\varepsilon} = 0 $ for all $\varepsilon$, and the law of the process
is symmetric, which means that $t_0(\varepsilon,\eta) =
t_1(\varepsilon,\eta) = +\infty$ for all $\varepsilon>0$ and
$\eta>0$. Moreover, $\sigma^2_\varepsilon= 2 c \varepsilon$ and Lemma
\ref{Lmm:KTPBnd} implies that $\bbp[\sup_{s\leq t} X_t \geq\eta] \leq
t^{{\eta}/{\varepsilon}} C(\eta,\varepsilon)$ and $\bbp[\sup_{s\leq
t} |X_t| \geq\eta] \leq
2 t^{{\eta}/{\varepsilon}} C(\eta,\varepsilon)$ with
$C(\eta,\varepsilon) =  (\frac{2c\mathrm{e}}{\eta} )^{{\eta
}/{\varepsilon}}$. The results of the above section
can {now} be improved to
\begin{eqnarray*}
\mathrm{e}^{(0)}(0,y,t) &\leq&\|\varphi\|_\infty\frac{4 \mathrm{e}^{-\lambda_\varepsilon t}}{b-y} C
\bigl(\varepsilon,(b-y)/2\bigr) t^{{(b-y)}/{(2\varepsilon)}},
\\
\mathrm{e}^{(1,1)}(0,y,t) &\leq&\|\varphi\|_\infty \mathrm{e}^{-\lambda_\varepsilon t}
a_\varepsilon C(\varepsilon,b) t^{1+{b}/{\varepsilon}},
\\
\mathrm{e}^{(1,2)}(0,y,t) & \leq&2 \|\varphi\|_\infty \mathrm{e}^{-\lambda_\varepsilon t}
a_\varepsilon C\bigl(\varepsilon,(b-y)/2\bigr) t^{1+{(b-y)}/{(2\varepsilon)}},
\\
\mathrm{e}^{(2,1)}(0,y,t) & \leq&\frac{\|\varphi\|_\infty}{2}\mathrm{e}^{-\lambda
_\varepsilon t}
a_\varepsilon\lambda_\varepsilon C(b,\varepsilon)
t^{2+{b}/{\varepsilon}},
\\
\mathrm{e}^{(2,3)}(0,y,t) & \leq&\|\varphi\|_\infty \mathrm{e}^{-\lambda_\varepsilon t}
a_\varepsilon\lambda_\varepsilon C\bigl((b-y)/2,\varepsilon\bigr)
t^{2+{(b-y)}/{(2\varepsilon)}}.
\end{eqnarray*}
To estimate $\mathrm{e}^{(2,2)}$ more precisely, let $\varepsilon_0 < \frac
{b-y-\varepsilon}{2} \wedge
(b-\varepsilon)$. Then
\begin{eqnarray*}
\bigl|D^1_t(y)\bigr| &\leq&2\|\varphi\|_\infty
a_\varepsilon\lambda_\varepsilon \mathbb P \bigl(U^\varepsilon_t
\geq\varepsilon_0 \bigr)+ {2\|\varphi\|_{\mathrm
{Lip}} \mathbb E
\biggl[ U^\varepsilon_t\int_b^\infty
s_\varepsilon\bigl(v-\bar{X}^\varepsilon_t
\bigr)s_\varepsilon\bigl(y-v +\bar{X}^\varepsilon_t -
X^\varepsilon_t \bigr) \,\mathrm{d}v \biggr]}
\\
&&{} + \|\varphi\|_\infty\mathbb E \biggl[ \mathbf{\mathbf{1}}_{U^\varepsilon_t
< \varepsilon_0}
\int_b^\infty \bigl(s_\varepsilon\bigl(v-
\bar{X}^\varepsilon_t\bigr)s_\varepsilon\bigl(y-v +
\bar{X}^\varepsilon_t - X^\varepsilon_t \bigr)
- s_\varepsilon(v) s_\varepsilon(y-v)\bigr)\,\mathrm{d}v \biggr]
\\
&\leq&2a_\varepsilon\lambda_\varepsilon \bigl( \|\varphi
\|_\infty \mathbb P \bigl(U^\varepsilon_t \geq
\varepsilon_0 \bigr) + \|\varphi\|_{\mathrm{Lip}}\bbe
\bigl[U^\varepsilon_t\bigr] \bigr)
\\
&&{} - \|\varphi\|_\infty\mathbb E \bigl[ U^\varepsilon_t
\bigr]\int_b^\infty s'_\varepsilon(v-
\varepsilon_0)s_\varepsilon(y-v +2\varepsilon_0)\,\mathrm{d}v
\\
&&{} + 2\|\varphi\|_\infty\mathbb E \bigl[ U^\varepsilon_t
\bigr]\int_b^\infty s_\varepsilon(v)s'_\varepsilon
\biggl(y-v +\frac{b-y}{2} \biggr)\,\mathrm{d}v
\\
&\leq& 2a_\varepsilon\lambda_\varepsilon \bigl( \|\varphi
\|_\infty \mathbb P \bigl(U^\varepsilon_t \geq
\varepsilon_0 \bigr) + \|\varphi\|_{\mathrm{Lip}}\bbe
\bigl[U^\varepsilon_t\bigr] \bigr)\\
&&{} + 2\|\varphi
\|_\infty\mathbb E \bigl[ U^\varepsilon_t \bigr]
s_\varepsilon (b-\varepsilon_0 ) s_\varepsilon (b-y-2
\varepsilon_0 ).
\end{eqnarray*}
A similar argument shows that
\[
\bigl|D^2_t(y)\bigr| \leq s_\varepsilon(b)
\lambda_\varepsilon\bigl(2\|\varphi\|_{\mathrm{Lip}} \bbe
\bigl[U^\varepsilon_t\bigr] + \|\varphi\|_\infty\bbp
\bigl[\bar X^\varepsilon_t \geq b\bigr]\bigr) + \|\varphi\|
_\infty\mathbb E \bigl[ U^\varepsilon_t \bigr]
s_\varepsilon (b ) s_\varepsilon ({b-y} ),
\]
which means that the bound for $|D^1_t(y)|$ always dominates.
Using the former bound, {we finally find the following upper bound for
$\mathrm{e}^{(2,2)}(0,y,t)$:}
\begin{eqnarray*}
&&2 \|\varphi\|_\infty \mathrm{e}^{-\lambda_\varepsilon t} a_\varepsilon
\lambda_\varepsilon C(\varepsilon_0,\varepsilon)t^{2+{\varepsilon_0}/{\varepsilon}}
\\
&&\quad{}+2 \mathrm{e}^{-\lambda_\varepsilon t} t^{{5}/{2}} \sigma_\varepsilon \bigl
\{s_\varepsilon (b-\varepsilon_0 ) s_\varepsilon ({b-y - 2
\varepsilon_0} )\|\varphi\|_\infty+ a_\varepsilon
\lambda_\varepsilon\|\varphi\|_{\mathrm{Lip}}\bigr\}.
\end{eqnarray*}
To specialize the estimate $\mathrm{e}^{(3)}$, we {upper bound} $\lambda
_\varepsilon^n \mathbb P(\bar X_t \geq b, X_t \in I_\delta
(y) | N^\varepsilon_t = n) $ by
\begin{eqnarray*}
&&\lambda_\varepsilon^n \mathbb P \Biggl(\bar
X^\varepsilon_t + \max_{0\leq k \leq n} \sum
_{i=1}^k Y_i \geq b,
X^\varepsilon_t + \sum_{i=1}^n
Y_i \in I_\delta(t) \Biggr)
\\
&&\quad\leq\lambda_\varepsilon^n \sum_{k=0}^n
\mathbb P \Biggl(\bar X^\varepsilon_t + \sum
_{i=1}^k Y_i \geq b,
X^\varepsilon_t + \sum_{i=1}^n
Y_i \in I_\delta(t) \Biggr).
\end{eqnarray*}
The cases $k=0$ and $k=n$ are treated separately:
\begin{eqnarray*}
&&\lambda_\varepsilon^n \mathbb P \Biggl(\bar
X^\varepsilon_t \geq b, X^\varepsilon_t +
\sum_{i=1}^n Y_i \in
I_\delta(t) \Biggr)
\\
&&\quad \leq\delta\mathbb P \Bigl(\sup_{u\leq
t} X^\varepsilon_u
\geq b \Bigr) \sup_{x} s^{*n}_\varepsilon(x)
\leq\delta C(b,\varepsilon) t^{{b}/{\varepsilon}} \sup_{x}
s^{*n}_\varepsilon(x),
\\
&&\lambda_\varepsilon^n \mathbb P \Biggl(\bar
X^\varepsilon_t + \sum_{i=1}^n
Y_i\geq b, X^\varepsilon_t + \sum
_{i=1}^n Y_i\in I_\delta(t)
\Biggr)
\\
&&\quad \leq\delta\mathbb P \bigl( \bar X^\varepsilon_t -
X^\varepsilon_t + y + \delta\geq b \bigr) \sup
_{x} s^{*n}_\varepsilon(x)\leq2\delta C
\bigl((b-y)/2,\varepsilon\bigr) t^{{(b-y)}/{(2\varepsilon)}} \sup_{x}
s^{*n}_\varepsilon(x).
\end{eqnarray*}
For $0<k<n$,
\begin{eqnarray*}
&&\mathbb P \Biggl(\bar X^\varepsilon_t + \sum
_{i=1}^k Y_i \geq b,
X^\varepsilon_t + \sum_{i=1}^n
Y_i \in I_\delta(t) \Biggr)
\\
&&\quad= \mathbb E \biggl[\int_y^{y+\delta}\,\mathrm{d}u \int
_{b-\bar X^\varepsilon_t}^\infty \,\mathrm{d}v\, s^{*k}_\varepsilon(v)s^{*(n-k)}_\varepsilon
\bigl(u-v-X^\varepsilon_t\bigr) \biggr]
\\
&&\quad \leq\delta\sup_{x} s^{*n}_\varepsilon(x)
\mathbb P \bigl(U^\varepsilon_t \geq\varepsilon_0
\bigr) + \delta \bar s^{*k}_\varepsilon(b-\varepsilon_0)
\int_{b}^\infty \,\mathrm{d}v\, \bar s^{*(n-k)}_\varepsilon
(y-v+2\varepsilon_0 +\delta),
\end{eqnarray*}
where $\bar s_\varepsilon$ is any function which is increasing on
$(-\infty,0)$, decreasing on $(0,\infty)$ and satisfies $\bar
s_\varepsilon(x) \geq s_\varepsilon(x)$ for all $x$. For the Cauchy
process, one can take $\bar s_\varepsilon(x) = \frac{2c}{x^2 +
\varepsilon^2}$
so that
\[
\bar s^{*k}_\varepsilon(x) = \frac{1}{\uppi} \biggl(
\frac{2 \uppi c}{\varepsilon} \biggr)^k \frac{\varepsilon
k}{(\varepsilon k)^2 + x^2}, \qquad\int
_b ^\infty\bar s_\varepsilon^{*k}(v)\,\mathrm{d}v
= \frac{1}{\uppi} \biggl(\frac{2 \uppi c}{\varepsilon} \biggr)^k \arctan
\frac{\varepsilon k}{b}.
\]
{Assembling all the estimates together, we finally get
\begin{eqnarray*}
\mathrm{e}^{(3)}(0,y,t)&\leq& \frac{\|\varphi\|_\infty}{6} a_\varepsilon
\lambda_\varepsilon^2 t^{3} \bigl(C(b,\varepsilon)
t^{{b}/{\varepsilon}}
\\
&&{} + 2C\bigl((b-y)/2,\varepsilon\bigr)t^{{(b-y)}/{(2\varepsilon)}} + 2C(
\varepsilon_0,\varepsilon) t^{{\varepsilon_0}/{\varepsilon}}\bigr)
\\
&&{}+\frac{16 \uppi c^3 t^3 \|\varphi\|_\infty}{3 \varepsilon
(b-\varepsilon_0)^2
(b-y-2\varepsilon_0)} \mathrm{e}^{{2\uppi c t}/{\varepsilon} -
\lambda_\varepsilon t}.
\end{eqnarray*}
The above estimates satisfy condition \eqref{Eq:RmdBhv} for
$\varepsilon< \frac{b-y}{4}\wedge b$. In the numerical
examples discussed in the paper, we have taken $\varepsilon=
\frac{b-y}{8}\wedge\frac{b}{2}$ and $\varepsilon_0 =
\frac{b-y}{4}\wedge\frac{b}{2}$.}
\end{appendix}

\section*{Acknowledgements}
The authors are grateful to an anonymous referee and both the Associate
Editor and the
editor-in-chief for their constructive and insightful comments that
greatly helped to improve the paper.
Jos\'e E. Figueroa-L\'opez research was partially supported by grants
from the US National Science Foundation (DMS-0906919, DMS-1149692).



\printhistory

\end{document}